\documentclass[final,hidelinks,onefignum,onetabnum]{siamart220329}

\usepackage{lipsum}
\usepackage{amsfonts}
\usepackage{graphicx}
\usepackage{epstopdf}
\usepackage{algorithmic}
  \DeclareGraphicsExtensions{.eps}

\newcommand{\del}[1]{\Delta #1}
\newcommand{\dstyle}[1]{\displaystyle{#1}}
\newcommand{\bmath}[1]{\mbox{\boldmath$ #1 $}}
\newcommand{\defeq}{\stackrel{\rm def}{=}}
\renewcommand{\d}{{\rm d}}
\newcommand{\e}{{\rm e}}
\newsiamremark{remark}{Remark}
\newsiamremark{hypothesis}{Hypothesis}
\crefname{hypothesis}{Hypothesis}{Hypotheses}
\newsiamthm{claim}{Claim}

\headers{Formulae for mixed moments
}{Y. Komori, G. Yang, and K. Burrage}

\title{Formulae for mixed moments of Wiener processes
  and a stochastic area integral
  \thanks{Submitted to the editors September 5, 2022.
    \funding{This work was partially supported by JSPS Grant-in-Aid 
      for Scientific Research~22K03416.}
  }
}

\author{Yoshio Komori\thanks{Department of Physics and Information
    Technology, Kyushu Institute of Technology, Japan
  (\email{komori@phys.kyutech.ac.jp}, \email{yoshio.kom@gmail.com}).}
\and Guoguo Yang\thanks{LMAM and School of Mathematical Sciences,
  Peking University, China
  (
  \email{yangguoguo@math.pku.edu.cn}}, \email{ygj512@hotmail.com}).
\and Kevin Burrage\thanks{School of Mathematical Sciences,
  Queensland University of Technology, Australia
  (\email{kevin.burrage@qut.edu.au}, \email{kevin.burrage@gmail.com}).}
}

\usepackage{amsopn}

\begin{document}

\maketitle

\begin{abstract}
This paper deals with the expectation of monomials
with respect to the stochastic area integral
\[
A_{1,2}(t,t+h)=\int_{t}^{t+h}\int_{t}^{s}\d W_{1}(r)\d W_{2}(s)
-\int_{t}^{t+h}\int_{t}^{s}\d W_{2}(r)\d W_{1}(s)
\]
and the increments of two Wiener processes,
$\Delta{W}_{i}(t,t+h)=W_{i}(t+h)-W_{i}(t),\ i=1,2$.
In a monomial, if the exponent of one of the Wiener increments
or the stochastic area integral is an odd number,
then the expectation of the monomial is zero.
However, if the exponent of any of them is an even number, then
the expectation is nonzero and its exact value is not known
in general.
In the present paper, we derive formulae to give the
value in general.
As an application of the formulae, we will utilize
the formulae for a careful stability analysis
on a Magnus-type Milstein method.
As another application, we will give
some mixed moments of the increments of Wiener processes
and stochastic double integrals.
\end{abstract}



\section{Introduction}
\label{sec:intro}
We are concerned with developing and analyzing numerical
methods that give strong first order approximations to
the solution of noncommutative stochastic differential equations (SDEs).
Such methods are usually constructed on the basis of
the comparison with the It\^{o}--Taylor expansion or the Stratonovich--Taylor
expansion, and as a result, they have one or more terms
related to the stochastic area integral
\[
A_{i,j}(t,t+h)=\int_{t}^{t+h}\int_{t}^{s}\d W_{i}(r)\d W_{j}(s)
-\int_{t}^{t+h}\int_{t}^{s}\d W_{j}(r)\d W_{i}(s),
\]
where $t\geq 0$, $h>0$,
and where $W_{i}(t),W_{j}(t)$ are independent Wiener processes
for positive integers $i,j\ (i\neq j)$.

L\'{e}vy \cite{Levy:1951} has studied the stochastic area integral
by utilizing the Fourier series of Wiener process.
As an example, he has given the probability density function
of $A_{i,j}(0,1)$. Gaveau \cite{Gaveau:1977} has also
studied the stochastic area integral and shown a joint
density function in general form, which is related
to Wiener increments and stochastic area integrals.

The joint density of the stochastic area integral $A_{1,2}(t,t+h)$
and the Wiener increments
$\Delta{W}_{i}(t,t+h)=W_{i}(t+h)-W_{i}(t),\ i=1,2$,
can be used to generate random numbers which
approximate the stochastic area integral.
In fact, by utilizing a conditional joint density function
of them, Gaines and Lyons \cite{Gaines:1994}
have proposed a random generator
for the stochastic area integral.
As the performance of random generators for the stochastic area integral
is a very important issue in stochastic simulation, many researchers have
struggled to find efficient generators
\cite{Kloeden:1992a,Kuznetsov:2018,Kuznetsov:2019,Mrongowius:2021,Wiktorsson:2001}.

Let us suppose that we have to calculate the expectation of a function
$f(y_{n+1})$ to study properties of a numerical method for
approximating the solution of SDEs.
Here, $y_{n+1}$ denotes an approximate solution
at time $t=t_{n+1}=t_{n}+h$.
If the numerical method contains
$I_{i,j}(t_{n},t_{n+1})=(1/2)\Delta{W}_{i}(t_{n},t_{n+1})\Delta{W}_{j}(t_{n},t_{n+1})
+(1/2)A_{i,j}(t_{n},t_{n+1})$ and
if we use the Taylor expansion of $f(y_{n+1})$ centered at $y_{n}$,
then we will need to calculate the expectation of monomials
with respect to $\Delta{W}_{i}(t_{n},t_{n+1})$,
$\Delta{W}_{j}(t_{n},t_{n+1})$, and $A_{i,j}(t_{n},t_{n+1})$.
Kloeden and Platen \cite{Kloeden:1999} have given
the expectation of monomials
with respect to $\Delta{W}_{i}(t_{n},t_{n+1})$,
$\Delta{W}_{j}(t_{n},t_{n+1})$, and $I_{i,j}(t_{n},t_{n+1})$
in some restricted cases.
However, to the best of the authors' knowledge,
the exact value of such expectations is not given in general.
Thus, in the present paper we will aim to derive
a formula that gives the exact value in general.

Incidentally, Magnus-type expansions for SDEs have recently drawn attention
from a number of researchers. Originally, a Magnus expansion
has been given for ordinary differential equations
by Magnus \cite{Magnus:1954}.
However, the Magnus expansion has been extended for SDEs, and
its extended expansions
have been studied by many researchers
\cite{Burrage:1999,Kamm:2020,Tambue:2020,Wang:2020,Wang:2020a}.
Especially, Yang et al. \cite{Yang:2021} have proposed
a Magnus-type Euler method and a Magnus-type Milstein method
for semilinear noncommutative It\^{o} SDEs.
When we investigate stability properties of the Magnus-type
Milstein method for a noncommutative test SDE,
we need the expectation of monomials
with respect to $\Delta{W}_{i}(t_{n},t_{n+1})$,
$\Delta{W}_{j}(t_{n},t_{n+1})$, and $A_{i,j}(t_{n},t_{n+1})$
in general.
Thus, as an application of the above formula, we will
carefully analyse stability properties of the method
for a noncommutative test SDE.

The present paper is organized as follows.
In \cref{sec:main_theorem}, after preliminary discussion we will
give our main theorem.
In \cref{sec:application}, we will introduce
the Magnus-type methods and a noncommutative test equation,
and as an application of our formula,
we will show that the Magnus-type Milstein method cannot be A-stable,
whereas the Magnus-type Euler method is A-stable for this test equation.
Finally, we will give concluding remarks.

\section{Main theorem}
\label{sec:main_theorem}
\subsection{Preliminary and motivation}
\label{subsec:preliminary}

For simplicity, let us take variables $a$,
$w_{i}$
corresponding to $A_{1,2}(t,t+1)$ and $\del{W}_{i}(t,t+1)$ $(i=1,2)$.
Then, the joint density function for $A_{1,2}(t,t+1)$, $\del{W}_{1}(t,t+1)$
and $\del{W}_{2}(t,t+1)$ is 
given as follows \cite{Gaines:1994}:
\begin{equation}
  f(a,
  w_{1},w_{2}
  )=\frac{1}{2\pi^{2}}\int_{0}^{\infty}
  \frac{x}{\sinh(x)}
  \exp\left(\frac{-(
    w_{1}^{2}+w_{2}^{2}
    )x}{2\tanh(x)}
  \right)\cos(ax)\d x.
  \label{eq:dens_a_bi}
\end{equation}
From this,
\begin{eqnarray}
  \lefteqn{
    \int_{-\infty}^{\infty}\int_{-\infty}^{\infty}
    w_{1}^{2k_{1}}w_{2}^{2k_{2}}
    f(a,w_{1},w_{2})
    \d w_{1}\d w_{2}
  }
  \nonumber
  \\
  &&
  =
  \frac{1}{\pi}\int_{0}^{\infty}
  \{(2k_{1}-1)!!\}
  \{(2k_{2}-1)!!\}
  \left(\frac{\tanh(x)}{x}
  \right)^{k_{1}+k_{2}}
  \frac{\cos(ax)}{\cosh(x)}\d x
  \nonumber
  \\
  &&
  =
  \frac{1}{\pi}
  \frac{(2k_{1})!((2k_{2})!)}{(k_{1}!)(k_{2}!)2^{k_{1}+k_{2}}}
  \int_{0}^{\infty}
  \left(\frac{\tanh(x)}{x}
  \right)^{k_{1}+k_{2}}
  \frac{\cos(ax)}{\cosh(x)}\d x
  \label{eq:marginal_dist_2k1_2k2}
\end{eqnarray}
for any nonnegative integers $k_{1},k_{2}$,
where $(2k-1)!!$ denotes $(2k-1)(2k-3)\cdots 1$ for a positive
integer $k$ and $(-1)!!=1$.
Thus, if we denote by $\gamma_{n,k,l}$ the expectation of monomial
$E[\{\del{W}_{1}(t,t+1)\}^{2n}
  \{A_{1,2}(t,t+1)\}^{2k}
  \{\del{W}_{2}(t,t+1)\}^{2l}
]$ for nonnegative integers $n,k$ and $l$,
then it can be expressed as
\begin{equation}
  \gamma_{n,k,l}
  =\frac{2}{\pi}
  \frac{((2n)!)((2l)!)}{(n!)(l!)2^{n+l}}
  \int_{0}^{\infty}
  \int_{0}^{\infty}
  a^{2k}
  \left(\frac{\tanh(x)}{x}
  \right)^{n+l}
  \frac{\cos(a x)}{\cosh(x)}\d x\d a.
  \label{eq:expect_w1w2DelI21_for_numer}
\end{equation}
From this, we have
\begin{equation}
  \gamma_{n-2k-l,k,l}
  =\frac{{}_{n-2k}C_{l}}{{}_{2(n-2k)}C_{2l}}
  \gamma_{n-2k,k,0}
  \label{eq:expect_w1w2DelI21_for_simp}
\end{equation}
for nonnegative integers $k,l$ and $n\geq 2k+l$,
where ${}_{n}C_{l}$ denotes a binomial coefficient
defined by $n!/((n-l)!(l!))$.
When $k=0$, we immediately have
\[
\gamma_{n-l,0,l}
=
\frac{(2(n-l))!((2l)!)}{(n-l)!(l!)2^{n}},
\]
which leads to
\[
\int_{0}^{\infty}
\int_{0}^{\infty}
\left(\frac{\tanh(x)}{x}
\right)^{n}
\frac{\cos(a x)}{\cosh(x)}\d x\d a
=\frac{\pi}{2}.
\]

When we want to consider $\gamma_{n-2k-l,k,l}$
for $k\neq 0$, we can concentrate on $\gamma_{n-2k,k,0}$
since this gives $\gamma_{n-2k-l,k,l}$ by
\eqref{eq:expect_w1w2DelI21_for_simp}, but
it is still not trivial to calculate it.
One possible way is to utilize the properties of stochastic integrals.
Let us denote by $I_{i,j}(t,t+1)$ the double stochastic integral
$\int_{t}^{t+1}\int_{t}^{s}\d W_{i}(r)\d W_{j}(s)$.
Noting that
\[
I_{1,2}(t,t+1)
=\frac{1}{2}\del{W}_{1}(t,t+1)\del{W}_{2}(t,t+1)+\frac{1}{2}A_{1,2}(t,t+1)
\]
(see \cite{Wiktorsson:2001}) and $A_{2,1}(t,t+1)=-A_{1,2}(t,t+1)$, we have
\begin{eqnarray*}
  &&
  I_{1,2}(t,t+1)+I_{2,1}(t,t+1)=\del{W}_{1}(t,t+1)\del{W}_{2}(t,t+1),
  \\
  &&
  I_{1,2}(t,t+1)-I_{2,1}(t,t+1)=A_{1,2}(t,t+1).
\end{eqnarray*}
From these,
\[
\bigl(A_{1,2}(t,t+1)\bigr)^{2}
=\bigl(\del{W}_{1}(t,t+1)\bigr)^{2}\bigl(\del{W}_{2}(t,t+1)\bigr)^{2}
-4I_{1,2}(t,t+1)I_{2,1}(t,t+1).
\]

Utilizing this and
\begin{equation}
  E\left[
    \bigl(\del{W}_{1}(t,t+1)\bigr)^{2}
    I_{1,2}(t,t+1)I_{2,1}(t,t+1)
    \right]
  =\frac{1}{3}
  \label{eq:ex_w1_integ}
\end{equation}
\cite[p. 225]{Kloeden:1999}, we can obtain
$\gamma_{1,1,0}=5/3$, which also gives $\gamma_{0,1,1}=5/3$
by \eqref{eq:expect_w1w2DelI21_for_simp}
for $n=3$ and $k=l=1$.

In the above calculations, \eqref{eq:ex_w1_integ} is a key point.
Such expectations are obtained from the properties of stochastic
integrals, but it is not easy to seek them when an exponent is a large
number. For example, in order to carry out calculations even for
\[
E\left[
  \bigl(\del{W}_{1}(t,t+1)\bigr)^{4}
  I_{1,2}(t,t+1)I_{2,1}(t,t+1)
  \right]
=2,
\]
we have needed a program code in a symbolic computing package, Mathematica,
which utilizes some rules for Stratonovich integrals $J_{12},J_{21}$
in \cite[p. 160]{Komori:2007}.
This fact motivates us to seek for another approach
in the next subsection.

\begin{remark}
  In contrast to the calculation of $\gamma_{n-2k-l,k,l}$, it is easy to get
  \[
  E\left[\{\del{W}_{1}(t,t+1)\}^{n}
    \{A_{1,2}(t,t+1)\}^{k}
    \{\del{W}_{2}(t,t+1\})^{l}
    \right]
  =0
  \]
  if one of $n$, $k$ and $l$ is an odd number.
  In fact, it is clear that
  \[
  \int_{-\infty}^{\infty}
  w_{1}^{2k_{1}+1}f(a,w_{1},w_{2})\d w_{1}
  =
  \int_{-\infty}^{\infty}
  w_{2}^{2k_{2}+1}f(a,w_{1},w_{2})\d w_{2}
  =0
  \]
  from \eqref{eq:dens_a_bi} and
  \[
  \int_{-\infty}^{\infty}\int_{-\infty}^{\infty}\int_{-\infty}^{\infty}
  w_{1}^{2k_{1}}w_{2}^{2k_{2}}a^{2k_{3}+1}
  f(a,w_{1},w_{2})
  \d w_{1}\d w_{2}
  \d a=0
  \]
  from \eqref{eq:marginal_dist_2k1_2k2}
  for any nonnegative integers $k_{1},k_{2},k_{3}$.
\end{remark}
%
%
\subsection{Lemmas and main theorem}
\label{subsec:theorem}

As preparation for our main theorem, we begin with a function
$r_{n}(x)$ given as
\[
r_{n}(x)=\left(\frac{\tanh(x)}{x}
\right)^{n}\frac{1}{\cosh(x)}
\]
for $n=0,1,\ldots$ and any nonzero real number $x$.
Recalling the series expansion
\[
\tanh(x)=x-\frac{x^{3}}{3}+\frac{2x^{5}}{15}
+\cdots
+\frac{2^{2n}(2^{2n}-1)B_{2n}}{(2n)!}x^{2n-1}+\cdots
\quad(x^{2}<\frac{\pi^{2}}{4})
\]
\cite[p. 42]{Zwillinger:2014},
where $B_{2n}$ is the $2n$-th Bernoulli number that
is given by
\[
B_{2n}=-\frac{1}{2n+1}+\frac{1}{2}
-\sum_{k=1}^{n-1}
\frac{2n(2n-1)\cdots(2n-2k+2)}
{(2k)!}B_{2k}\quad(n=1,2,\ldots)
\]
\cite[p. 1052]{Zwillinger:2014},
we define the value of $r_{n}(0)$ by
$\lim_{x\to 0}r_{n}(x)=1$.
If we differentiate $r_{n}(x)$, then
\begin{equation}
  \frac{\d r_{n}}{\d x}(x)=v_{n}(x)r_{n}(x),
  \label{eq:dr/dx}
\end{equation}
where
\[
v_{n}(x)=n\left(
\frac{1}{\tanh(x)}-\frac{1}{x}
\right)-(n+1)\tanh(x).
\]
Recalling the series expansion
\[
\frac{1}{\tanh(x)}
=\frac{1}{x}+\frac{x}{3}-\frac{x^{3}}{45}
+\frac{2x^{5}}{945}
+\cdots
+\frac{2^{2n}B_{2n}}{(2n)!}x^{2n-1}+\cdots
\quad(x^{2}<\pi^{2})
\]
\cite[p. 42]{Zwillinger:2014},
we define the value of $v_{n}(0)$ by $\lim_{x\to 0} v_{n}(x)=0$.
Similarly, we define the value
of $(\d^{2k} v_{n}/\d x^{2k})(0)$ by
$
\lim_{x\to 0} (\d^{2k} v_{n}/\d x^{2k})(x)=0
$
for $k=1,2,\ldots$.
Now, we can obtain the following lemma regarding $r_{n}(x)$.

\begin{lemma}
  \label{lem:rec_formula_rn}
  For $k=1,2,\ldots$,
  the derivatives of $r_{n}(x)$ at $x=0$ are given by
  \begin{eqnarray}
    &&
    \frac{\d^{2k-1}r_{n}}{\d x^{2k-1}}(0)=0,
    \nonumber
    \\
    &&
    \frac{\d^{2k}r_{n}}{\d x^{2k}}(0)
    =\frac{1}{2k}\sum_{j=1}^{k}{}_{2k}C_{2j}2^{2j}B_{2j}
    \left\{(1-2^{2j})(n+1)+n
    \right\}\frac{\d^{2(k-j)}r_{n}}{\d x^{2(k-j)}}(0).
    \label{eq:rec_formula_rn}
  \end{eqnarray}
\end{lemma}
\begin{proof}
As $v_{n}(0)=0$ and $r_{n}(0)=1$, \eqref{eq:dr/dx} gives
$(\d r_{n}/\d x)(0)=0$.
Next, by applying Leibniz theorem to \eqref{eq:dr/dx} and
noting that $(\d^{2k}v_{n}/\d x^{2k})(0)=0$ for
$k=1,2,\ldots$, we have
\[
\frac{\d^{2k-1}r_{n}}{\d x^{2k-1}}(0)=0,
\qquad
\frac{\d^{2k}r_{n}}{\d x^{2k}}(0)
=\sum_{j=1}^{k}{}_{2k-1}C_{2j-1}
\frac{\d^{2j-1}v_{n}}{\d x^{2j-1}}(0)
\frac{\d^{2(k-j)}r_{n}}{\d x^{2(k-j)}}(0).
\]
Here, as
\[
\frac{\d^{2j-1}v_{n}}{\d x^{2j-1}}(0)
=\frac{2^{2j}B_{2j}}{2j}
\left\{(1-2^{2j})(n+1)+n
\right\},
\]
the substitution of this and simplification yield
the second equation in the lemma.
\end{proof}

\begin{remark}
  \label{rem:limit_rn}
  If we set $y=1/\tanh(x)$, then $\d y/\d x=1-y^{2}$.
  Noting $\d /\d x=(1-y^{2})(\d /\d y)$  and $\lim_{x\to\infty}y=1$,
  we have $\lim_{x\to\infty}(\d^{k}y/\d x^{k})=0$
  for any positive integer $k$. We have the same equation
  also for $y=\tanh(x)$.
  From these, we have $\lim_{x\to\infty}(\d^{k}v_{n}/\d x^{k})(x)=0$.
  Thus, in a similar way to the proof of \cref{lem:rec_formula_rn}
  we can see that $(\d^{k}r_{n}/\d x^{k}(x)$ exponentially converges
  to $0$ as $x\to\infty$.
\end{remark}

As \cref{lem:rec_formula_rn} gives a recursive formula
for the even-order derivatives
of $r_{n}(x)$ at $x=0$, it is useful for the fast calculations of them.
On the other hand, the following lemma gives an explicit formula,
which can be helpful for mathematical analysis.

\begin{lemma}
  \label{lem:explicit_formula_rn}
  The even order derivatives of $r_{n}(x)$ at $x=0$ are given by
  \begin{equation}
    \frac{\d^{2k}r_{n}}{\d x^{2k}}(0)
    =(-1)^{k}\bigl((2k)!\bigr)s_{n,k}
    \label{eq:explicit_formula_rn}
  \end{equation}
  for $k=0,1,\ldots$. Here, $s_{n,k}$ is given by
  \[
  s_{n,0}=1,\qquad
  s_{n,k}
  =\sum_{l_{1}+2l_{2}+\cdots+kl_{k}=k}
  \left\{
  \prod_{j=1}^{k}\frac{\beta_{n,j}^{l_{j}}}{j^{l_{j}}(l_{j}!)}
  \right\}
  \quad(k=1,2,\ldots),
  \]
  where $l_{1},l_{2},\ldots,l_{k}$
  denote nonnegative integers
  and where
  \[
  \beta_{n,j}
  =\frac{2^{2j-1}|B_{2j}|
  \left\{(2^{2j}-1)(n+1)-n
  \right\}}{(2j)!}.
  \]
\end{lemma}
\begin{proof}
The setting of $k=0$ is a special case and it is trivial
that \eqref{eq:explicit_formula_rn} holds.
In what follows, we shall prove that \eqref{eq:explicit_formula_rn} holds
by a mathematical induction.
When $k=1$, \eqref{eq:rec_formula_rn} implies
\[
\frac{\d^{2}r_{n}}{\d x^{2}}(0)
=2B_{2}
\left\{(1-2^{2})(n+1)+n
\right\}
=-2\beta_{n,1}.
\]
Thus, \eqref{eq:explicit_formula_rn} holds for $k=1$.
Next, suppose that \eqref{eq:explicit_formula_rn} holds
for any $k$ less than or equal to an even number $q=2m$.
From \eqref{eq:rec_formula_rn}, we have
\begin{eqnarray}
  &&
  \frac{\d^{2(q+1)}r_{n}}{\d x^{2(q+1)}}(0)
  =\frac{1}{q+1}\sum_{j=1}^{q+1}{}_{2(q+1)}C_{2j}
  \bigl((2j)!
  \bigr)(-1)^{j}\beta_{n,j}\frac{\d^{2(q+1-j)}r_{n}}{\d x^{2(q+1-j)}}(0)
  \nonumber
  \\
  &&
  \makebox[5.6em]{}
  =(-1)^{q+1}\frac{(2(q+1))!}{q+1}(H_{1}+H_{2}),
  \label{eq:all_sum_for_rec_formula}
\end{eqnarray}
where
\begin{eqnarray*}
  &&
  H_{1}=
  \sum_{\hat{l}_{1}+2l_{2}+\cdots+ql_{q}=q+1}
  \left[
    \hat{l}_{1}
    \frac{\beta_{n,1}^{\hat{l}_{1}}}{\hat{l}_{1}!}
    \prod_{\stackrel{\scriptstyle j=1}{j \neq 1}}^{q}
    \frac{\beta_{n,j}^{l_{j}}}{j^{l_{j}}(l_{j}!)}
    \right]
  \nonumber
  \\
  &&
  \makebox[2.5em]{}
  +\sum_{l_{1}+2\hat{l}_{2}+3l_{3}+\cdots+(q-1)l_{q-1}=q+1}
  \left[
    2\hat{l}_{2}
    \frac{\beta_{n,2}^{\hat{l}_{2}}}{2^{}\hat{l}_{2}(\hat{l}_{2}!)}
    \prod_{\stackrel{\scriptstyle j=1}{j \neq 2}}^{q-1}
    \frac{\beta_{n,j}^{l_{j}}}{j^{l_{j}}(l_{j}!)}
    \right]
  +\cdots
  \nonumber
  \\
  &&
  \makebox[2.5em]{}
  +\sum_{l_{1}+2l_{2}+\cdots+(m-1)l_{m-1}+m\hat{l}_{m}+(m+1)l_{m+1}=q+1}
  \left[
    m\hat{l}_{m}
    \frac{\beta_{n,m}^{\hat{l}_{m}}}{m^{\hat{l}_{m}}(\hat{l}_{m}!)}
    \prod_{\stackrel{\scriptstyle j=1}{j \neq m}}^{m+1}
    \frac{\beta_{n,j}^{l_{j}}}{j^{l_{j}}(l_{j}!)}
    \right],
  \nonumber
  \\
  &&
  \raisebox{0pt}[25pt][0pt]{}
  H_{2}=
  \sum_{l_{1}+2l_{2}+\cdots+ml_{m}+(m+1)\hat{l}_{m+1}=q+1}
  \left[
    \beta_{n,m+1}
    \prod_{j=1}^{m}\frac{\beta_{n,j}^{l_{j}}}{j^{l_{j}}(l_{j}!)}
    \right]
  \nonumber
  \\
  &&
  \makebox[2.5em]{}
  +\sum_{l_{1}+2l_{2}+\cdots+(m-1)l_{m-1}+(m+2)\hat{l}_{m+2}=q+1}
  \left[
    \beta_{n,m+2}
    \prod_{j=1}^{m-1}\frac{\beta_{n,j}^{l_{j}}}{j^{l_{j}}(l_{j}!)}
    \right]
  +\cdots
  \nonumber
  \\
  &&
  \makebox[2.5em]{}
  +\sum_{l_{1}+q\hat{l}_{q}=q+1}
  \left[
    \beta_{n,q}
    \prod_{j=1}^{1}\frac{\beta_{n,j}^{l_{j}}}{j^{l_{j}}(l_{j}!)}
    \right]
  +\beta_{n,q+1},
\end{eqnarray*}
and where $\hat{l}_{1},\hat{l}_{2},\ldots,\hat{l}_{q}$
denote positive integers.
In $H_{2}$, for example, 
we can rewrite the sum of the last two terms as
\begin{eqnarray*}
  &&
  \makebox[1em]{}
  \sum_{l_{1}+q\hat{l}_{q}=q+1}
  \left[
    \beta_{n,q}
    \prod_{j=1}^{1}\frac{\beta_{n,j}^{l_{j}}}{j^{l_{j}}(l_{j}!)}
    \right]
  +\beta_{n,q+1}
  \\
  &&
  =\sum_{l_{1}+ql_{q}+(q+1)l_{q+1}=q+1}
  \left\{
  \bigl(ql_{q}+(q+1)l_{q+1}
  \bigr)
  \left(
  \frac{\beta_{n,1}^{l_{1}}}{l_{1}!}
  \right)
  \left(
  \frac{\beta_{n,q}^{l_{q}}}{q^{l_{q}}(l_{q}!)}
  \right)
  \left(
  \frac{\beta_{n,q+1}^{l_{q+1}}}{(q+1)^{l_{q+1}}(l_{q+1}!)}
  \right)
  \right\}.
\end{eqnarray*}
Noting this, we can see that $H_{2}$ is rewritten as
\begin{eqnarray*}
  &&
  H_{2}=
  \sum_{l_{1}+2l_{2}+\cdots+(q+1)l_{q+1}=q+1}
  \left\{
  \raisebox{0pt}[20pt][0pt]{}
  \bigl((m+1)l_{m+1}+(m+2)l_{m+2}+\cdots
  \bigr.
  \right.
  \\
  &&
  \makebox[12.5em]{}
  \left.
  \bigl.
  +(q+1)l_{q+1}
  \bigr)\prod_{j=1}^{q+1}\frac{\beta_{n,j}^{l_{j}}}{j^{l_{j}}(l_{j}!)}
  \right\}.
\end{eqnarray*}
Similarly, $H_{1}$ is rewritten as
\begin{eqnarray*}
  &&
  H_{1}
  =\sum_{l_{1}+2l_{2}+\cdots+ql_{q}=q+1}
  \left\{
  \left(l_{1}+2l_{2}+\cdots+ml_{m}
  \right)\prod_{j=1}^{q}\frac{\beta_{n,j}^{l_{j}}}{j^{l_{j}}(l_{j}!)}
  \right\}
  \\
  &&
  \makebox[1em]{}
  =
  \sum_{l_{1}+2l_{2}+\cdots+(q+1)l_{q+1}=q+1}
  \left\{
  \left(l_{1}+2l_{2}+\cdots+ml_{m}
  \right)\prod_{j=1}^{q+1}\frac{\beta_{n,j}^{l_{j}}}{j^{l_{j}}(l_{j}!)}
  \right\}.
\end{eqnarray*}
The substitution of these into \eqref{eq:all_sum_for_rec_formula}
implies
\[
\frac{\d^{2(q+1)}r_{n}}{\d x^{2(q+1)}}(0)
=
(-1)^{q+1}\bigl((2(q+1))!\bigr)
\sum_{l_{1}+2l_{2}+\cdots+(q+1)l_{q+1}=q+1}
\left\{
\prod_{j=1}^{q+1}\frac{\beta_{n,j}^{l_{j}}}{j^{l_{j}}(l_{j}!)}
\right\}.
\]
Thus, \eqref{eq:explicit_formula_rn} holds for $k=q+1$.

On the other hand, when we suppose that \eqref{eq:explicit_formula_rn} holds
for any $k$ less than or equal to an odd number $q=2m+1$,
we can obtain the above equality for $k=q+1$ in a similar way.
%
\end{proof}

Rather than the explicit formula that gives the exact value of
even-order derivatives of $r_{n}(x)$ at $x=0$,
a simpler formula which gives
an upper
bound of the absolute value of them
can be often helpful for mathematical analysis.
As preparation, let us consider the summation
\[
\hat{s}_{n,L}(k)
=\sum_{l_{1}+2l_{2}+\cdots+Ll_{L}=k}
\left\{
\prod_{j=1}^{L}\frac{\beta_{n,j}^{l_{j}}}{j^{l_{j}}(l_{j}!)}
\right\}
\]
for nonnegative integers $L$ and $k$, which implies
$\hat{s}_{n,L}(0)=1$ and $\hat{s}_{n,0}(k)=0$ for $k>0$.
When we deal with $\hat{s}_{n,L}(k)$, the following lemma
is useful.

\begin{lemma}
  \label{lem:sum_sHat}
  Let $L$ be a nonnegative integer
  and $k$ a positive integer. Then,
  \[
  \hat{s}_{n,L}(k)
  =
  \left\{
  \begin{array}{ll}
    \raisebox{0pt}[22pt][22pt]{}
    \dstyle{
      \frac{1}{k}
      \sum_{j=1}^{L}\beta_{n,j}\hat{s}_{n,L}(k-j)
    }
    & \qquad (L<k),
    \\
    \dstyle{
      \frac{1}{k}
      \sum_{j=1}^{k}\beta_{n,j}\hat{s}_{n,k}(k-j)
    }
    & \qquad (L\geq k).
  \end{array}
  \right.
  \]
\end{lemma}
For the proof of this lemma, we refer the reader to
\cref{app:proof_lemma_sum_sHat}.
Using this, we obtain a function related
to $\hat{s}_{n,L}(k)$ as we can see below.

\begin{lemma}
  \label{lem:generating_func}
  For a nonnegative integer $L$, define $M_{n,L}(\theta)$ by
  \begin{equation}
    M_{n,L}(\theta)
    =\exp\left(\sum_{j=1}^{L}\frac{\beta_{n,j}}{j}\theta^{j}
    \right),
    \label{eq:generating_func}
  \end{equation}
  where $\theta\in\mathbb{R}$ is a parameter.
  For any nonnegative integer $k$,
  this satisfies
  \begin{equation}
  \hat{s}_{n,L}(k)=\frac{1}{k!}
  \frac{\d^{k}M_{n,L}}{\d\theta^{k}}(0).
  \label{eq:sHat_by_MnL}
  \end{equation}
\end{lemma}
\begin{proof}
When $L=0$ or $k=0$, it is trivial that \eqref{eq:sHat_by_MnL} holds.
When $L>0$ and $k>0$, we shall prove this lemma by a mathematical induction.
The differentiation of \eqref{eq:generating_func}
and the substitution of $\theta=0$ into it yield
$\frac{\d M_{n,L}}{\d \theta}(0)
=\beta_{n,1}=\hat{s}_{n,L}(1)$,
which implies that \eqref{eq:sHat_by_MnL} holds for $k=1$.
Next, by differentiating \eqref{eq:generating_func} and
applying Leibniz theorem to it,
we have
\begin{equation}
  \frac{\d^{k+1} M_{n,L}}{\d \theta^{k+1}}(\theta)
  =
  \sum_{l=0}^{k}{}_{k}C_{l}
  \left\{
  \frac{\d^{l}}{\d \theta^{l}}
  \left(\sum_{j=1}^{L}\beta_{n,j}\theta^{j-1}
  \right)
  \right\}
  \frac{\d^{k-l} M_{n,L}}{\d \theta^{k-l}}(\theta).
  \label{eq:k1_order_dMnL}
\end{equation}
Let us suppose that \eqref{eq:sHat_by_MnL} holds for any $k$
less than or equal to a positive integer $q$,
and consider the case of $q<L$.
Then, \eqref{eq:k1_order_dMnL} gives
\begin{eqnarray*}
  &&
  \frac{\d^{q+1} M_{n,L}}{\d \theta^{q+1}}(0)
  =
  \sum_{l=0}^{q}{}_{q}C_{l}
  \beta_{n,l+1}(l!)
  \frac{\d^{q-l} M_{n,L}}{\d \theta^{q-l}}(0)
  \\
  &&
  \makebox[5.5em]{}
  =
  \sum_{l=0}^{q}{}_{q}C_{l}
  \beta_{n,l+1}(l!)
  ((q-l)!)\hat{s}_{n,L}(q-l)
  \\
  &&
  \makebox[5.5em]{}
  =
  q!
  \sum_{l=1}^{q+1}
  \beta_{n,l}
  \hat{s}_{n,L}(q+1-l)
  \\
  &&
  \makebox[5.5em]{}
  =
  q!
  \sum_{l=1}^{q+1}
  \beta_{n,l}
  \hat{s}_{n,q+1}(q+1-l)
  =(q+1)!\hat{s}_{n,L}(q+1).
\end{eqnarray*}
In the last two equalities, we have used $L\geq q+1$ and
\cref{lem:sum_sHat}.
On the other hand, in the case of $q\geq L$,
\eqref{eq:k1_order_dMnL} similarly gives
\begin{eqnarray*}
  &&
  \frac{\d^{q+1} M_{n,L}}{\d \theta^{q+1}}(0)
  =
  \sum_{l=0}^{L-1}{}_{q}C_{l}
  \beta_{n,l+1}(l!)
  \frac{\d^{q-l} M_{n,L}}{\d \theta^{q-l}}(0)
  \\
  &&
  \makebox[5.5em]{}
  =
  q!
  \sum_{l=1}^{L}
  \beta_{n,l}
  \hat{s}_{n,L}(q+1-l)
  =(q+1)!\hat{s}_{n,L}(q+1).
\end{eqnarray*}
Thus, \eqref{eq:sHat_by_MnL} holds for $k=q+1$.
%
\end{proof}

Utilizing this lemma, the following lemma gives
upper bounds of $s_{n,k}$.

\begin{lemma}
  \label{lem:bounds_snk}
  Let $k_{0}$ be a nonnegative integer.
  Then,
  \begin{equation}
    s_{n,k}\leq
    \frac{\tan^{n}(1)}{\cos(1)}
    \frac{1}{M_{n,k_{0}}(1)}
    \max_{0\leq j\leq k}
    \left\{
    \frac{1}{j!}
    \frac{\d^{j}M_{n,k_{0}}}{\d\theta^{j}}(0)
    \right\}
    \label{ineq:bound_snk}
  \end{equation}
  for $k\geq k_{0}+1$.
  Especially, for $k\geq 0$,
  \begin{equation}
    s_{n,k}\leq
    \frac{\tan^{n}(1)}{\cos(1)}.
    \label{ineq:bound_simple_snk}
  \end{equation}
\end{lemma}
\begin{proof}
Setting $\hat{s}_{n,L}(k)=0$ if $k<0$, we have
\[
\hat{s}_{n,L}(k)=
\sum_{i=0}^{k}\left\{
\frac{\beta_{n,L}^{i}}{L^{i}(i!)}\hat{s}_{n,L-1}(k-Li)
\right\}
\]
for $L\geq 1$.
In a similar way, for $L\geq k_{0}+1$,
\begin{eqnarray*}
  &&
  \hat{s}_{n,L}(k)=
  \sum_{i_{1},i_{2},\ldots,i_{L-k_{0}}=0}^{k}\left\{
  \left(
  \prod_{j=1}^{L-k_{0}}
  \frac{\beta_{n,L-j+1}^{i_{j}}}{(L-j+1)^{i_{j}}(i_{j}!)}
  \right)
  \hat{s}_{n,k_{0}}\left(k-\sum_{j=1}^{L-k_{0}}(L-j+1)i_{j}
  \right)
  \right\}
  \\
  &&
  \makebox[3em]{}
  \leq
  \sum_{i_{1},i_{2},\ldots,i_{L-k_{0}}=0}^{k}
  \left(
  \prod_{j=1}^{L-k_{0}}
  \frac{\beta_{n,L-j+1}^{i_{j}}}{(L-j+1)^{i_{j}}(i_{j}!)}
  \right)
  \max_{0\leq j\leq k}
  \hat{s}_{n,k_{0}}(k-j)
  \\
  &&
  \makebox[3em]{}
  \leq
  \exp
  \left(
  \sum_{i=k_{0}+1}^{L}
  \frac{\beta_{n,i}}{i}
  \right)
  \max_{0\leq j\leq k}
  \hat{s}_{n,k_{0}}(j)
  \leq
  \frac{\tan^{n}(1)}{\cos(1)}
  \frac{1}{M_{n,k_{0}}(1)}
  \max_{0\leq j\leq k}
  \left\{
  \frac{1}{j!}
  \frac{\d^{j}M_{n,k_{0}}}{\d\theta^{j}}(0)
  \right\}.
\end{eqnarray*}
In the last inequality,
we have utilized \cref{lem:generating_func} and the formulae
\cite[p. 55]{Zwillinger:2014}
\[
\sum_{j=1}^{\infty}
\frac{2^{2j-1}(2^{2j}-1)|B_{2j}|}{j(2j)!}
=-\ln\cos(1),
\qquad
\sum_{j=1}^{\infty}
\frac{2^{2j-1}|B_{2j}|}{j(2j)!}
=-\ln\sin(1).
\]
Thus, when $L=k$ we obtain \eqref{ineq:bound_snk}.
In addition, the substitution of $k_{0}=0$ into \eqref{ineq:bound_snk}
and $s_{n,0}=1$ complete the proof.
\end{proof}

As the next step of preparations for our main theorem,
we shall derive
\begin{equation}
  \int_{0}^{\infty}\int_{0}^{\infty}r_{n}(x)\cos(ax)\d x\d a=\frac{\pi}{2}
  \label{eq:integ_rn_cos}
\end{equation}
in a different way, although it has already been obtained
in \cref{subsec:preliminary}.

When $n=0$, Levy \cite{Levy:1951} has shown
\[
\int_{0}^{\infty}r_{0}(x)\cos(ax)\d x=\frac{\pi}{2\cosh(\pi a/2)}
\]
by utilizing an expansion of $1/\cosh(\pi a/2)$
in series of simple fractions \cite[p. 44]{Zwillinger:2014}.
This directly gives \eqref{eq:integ_rn_cos} when $n=0$.
Now, let us consider another approach to calculate the double integral
when $n=0$.
Noting that we may interchange the integration and differentiation
with respect to $x$ and $a$, respectively in the following,
we obtain
\begin{eqnarray}
  &&
  \int_{0}^{\infty}\int_{0}^{\infty}r_{0}(x)\cos(ax)\d x\d a
  =\int_{0}^{\infty}\int_{0}^{\infty}r_{0}(x)
  \frac{\d}{\d a}
  \left(\frac{\sin(ax)}{x}
  \right)\d x\d a
  \nonumber
  \\
  &&
  \makebox[11.7em]{}
  =\lim_{a\to\infty}
  \int_{0}^{\infty}r_{0}(x)
  \frac{\sin(ax)}{x}
  \d x.
  \label{eq:interchange_r0}
\end{eqnarray}
Here, recalling one of the well-known results by the residue theorem:
\[
\int_{0}^{\infty}\frac{\sin(ax)}{x}\d x=\frac{\pi}{2},
\]
we consider
\[
\int_{0}^{\infty}r_{0}(x)
\frac{\sin(ax)}{x}\d x-\frac{\pi}{2}
=\int_{0}^{\infty}
\frac{r_{0}(x)-r_{0}(0)}{x}\sin(ax)\d x,
\]
which is interpreted as an improper Riemann integral.
Although we cannot directly apply the Riemann-Lebesgue
lemma \cite[p. 169]{Bachmann:2000} to the integral,
we can use its basic idea.
Taking \cref{lem:rec_formula_rn} and \cref{rem:limit_rn}
into account, we have
\[
a\int_{0}^{\infty}
\frac{r_{0}(x)-r_{0}(0)}{x}\sin(ax)\d x
=\int_{0}^{\infty}
\frac{r^{\prime}_{0}(x)x-r_{0}(x)+r_{0}(0)}{x^{2}}\cos(ax)\d x
<\infty.
\]
Thus, as
\[
\lim_{a\to\infty}
\left\{
\int_{0}^{\infty}r_{0}(x)
\frac{\sin(ax)}{x}\d x-\frac{\pi}{2}
\right\}=0,
\]
we obtain \eqref{eq:integ_rn_cos} when $n=0$.
Also for $n\geq 1$, we can show
\[
\lim_{a\to\infty}
\left\{
\int_{0}^{\infty}r_{n}(x)
\frac{\sin(ax)}{x}\d x-\frac{\pi}{2}
\right\}=0
\]
in a similar way. Thus, \eqref{eq:integ_rn_cos} holds
for $n=0,1,\ldots$.

We have already mentioned that
in order to obtain $\gamma_{n-2k-l,k,l}$,
we can concentrate on $\gamma_{n-2k,k,0}$
due to \eqref{eq:expect_w1w2DelI21_for_simp}.
Our main theorem is useful to get this expression.
\begin{theorem}
  \begin{equation}
    \gamma_{n,k,0}
    =(-1)^{k}
    \frac{(2n)!}{2^{n}(n!)}
    \frac{\d^{2k}r_{n}}{\d x^{2k}}(0)
    \label{eq:main_th_rec}
    =
    \frac{(2n)!((2k)!)}{2^{n}(n!)}
    s_{n,k}
    =
    \frac{(2n)!}{2^{n}(n!)}
    \frac{(2k)!}{k!}
    \frac{\d^{k}M_{n,k}}{\d\theta^{k}}(0)
    \end{equation}
  for any nonnegative integers $n,k$.
\end{theorem}
\begin{proof}
  The others except for the first equality
  are immediately given
by \cref{lem:explicit_formula_rn}
and \cref{lem:generating_func}
if the first equality holds,
all that remains is its proof.
For a positive integer $k$, we obtain
\[
\int_{0}^{\infty}\int_{0}^{\infty}a^{2k}r_{n}(x)\cos(ax)\d x\d a
=-\int_{0}^{\infty}\int_{0}^{\infty}
a^{2k-1}
\frac{\d r_{n}}{\d x}(x)
\sin(ax)\d x\d a,
\]
using integration by parts with respect to $x$
and noting \cref{rem:limit_rn}.
Its repeated applications and
the interchange of the integration
and differentiation lead to
\[
\int_{0}^{\infty}\int_{0}^{\infty}a^{2k}r_{n}(x)\cos(ax)\d x\d a
=\lim_{a\to\infty}(-1)^{k}
\int_{0}^{\infty}
\frac{\d^{2k} r_{n}}{\d x^{2k}}(x)
\frac{\sin(ax)}{x}
\d x
\]
in a similar way to \eqref{eq:interchange_r0}.
Here, taking \cref{lem:rec_formula_rn} and \cref{rem:limit_rn}
into account, we have
\begin{eqnarray*}
  &&
  \makebox[1em]{}
  a\int_{0}^{\infty}
  \frac{\frac{\d^{2k} r_{n}}{\d x^{2k}}(x)
    -\frac{\d^{2k} r_{n}}{\d x^{2k}}(0)}{x}\sin(ax)\d x
  \\
  &&
  =\int_{0}^{\infty}
  \frac{\frac{\d^{2k+1} r_{n}}{\d x^{2k}}(x)x
    -\frac{\d^{2k} r_{n}}{\d x^{2k}}(x)
    +\frac{\d^{2k} r_{n}}{\d x^{2k}}(0)}{x^{2}}\cos(ax)\d x
  <\infty.
\end{eqnarray*}
Thus,
\[
\lim_{a\to\infty}
\left\{
\int_{0}^{\infty}\frac{\d^{2k} r_{n}}{\d x^{2k}}(x)
\frac{\sin(ax)}{x}\d x-\frac{\pi}{2}\frac{\d^{2k} r_{n}}{\d x^{2k}}(0)
\right\}=0.
\]
Consequently,
\[
\int_{0}^{\infty}\int_{0}^{\infty}a^{2k}r_{n}(x)\cos(ax)\d x\d a
=(-1)^{k}
\frac{\pi}{2}\frac{\d^{2k} r_{n}}{\d x^{2k}}(0).
\]
The substitution of the equation and $l=0$ into
\eqref{eq:expect_w1w2DelI21_for_numer}
completes the proof.
\end{proof}

\begin{remark}
  In the interpretation of \eqref{eq:main_th_rec},
  $M_{n,k}(\theta)$ is the mixed moment generating function of
  $\{\del{W}_{1}(t,t+1)\}^{2}$
  and $\{A_{1,2}(t,t+1)\}^{2}$,
  whereas $M_{n,2}(\theta)$ is the moment generating function of
  a normal random variable $Y_{n}$ with
  mean $\beta_{n,1}$ and variance $\beta_{n,2}$
  in the usual meaning.
\end{remark}

This theorem and \cref{lem:bounds_snk} immediately
give the following corollary.
\begin{corollary}
  Let $k_{0}$ be a nonnegative integer.
  Then,
  \[
  \gamma_{n,k,0}\leq
  \frac{(2n)!((2k)!)}{2^{n}(n!)}
  \frac{\tan^{n}(1)}{\cos(1)}
  \frac{1}{M_{n,k_{0}}(1)}
  \max_{0\leq j\leq k}
  \left\{
  \frac{1}{j!}
  \frac{\d^{j}M_{n,k_{0}}}{\d\theta^{j}}(0)
  \right\}
  \]
  for $k\geq k_{0}+1$.
  Especially, for $k\geq 0$,
  \[
  \gamma_{n,k,0}\leq
  \frac{(2n)!((2k)!)}{2^{n}(n!)}
  \frac{\tan^{n}(1)}{\cos(1)}.
  \]
\end{corollary}

\section{Applications}
\label{sec:application}

We give two examples as applications of our main theorem.
The first example is the stability analysis for Magnus-type
methods. The second example is the mixed moments
of Wiener increments and stochastic double integrals.

\subsection{Stability analysis for Magnus-type methods}

Using a noncommutative test SDE,
we analyse stability properties of Magnus-type methods.
Then, a very careful treatment is necessary especially
for the Magnus-type Milstein method.
For this, first we derive polynomials that play
an important role in the analysis.

\subsubsection{Important polynomials in the analysis}

As a new type of method for SDEs, Magnus-type methods \cite{Yang:2021}
have been recently proposed for semilinear SDEs given by
\begin{equation}
  \d\bmath{y}(t)=
  \bigl\{F_{0}\bmath{y}(t)+\bmath{g}_{0}(\bmath{y}(t))
  \bigr\}\d t
  +\sum_{j=1}^{m}
  \bigl\{F_{j}\bmath{y}(t)+\bmath{g}_{j}(\bmath{y}(t))
  \bigr\}\d W_{j}(t),
  \quad \bmath{y}(0)=\bmath{y}_{0},
  \label{eq:sde}
\end{equation}
where $t\in[0,T]$ and where $\bmath{g}_{j}$,
$j=0,1,\ldots,m$ are
$\mathbb{R}^{d}$-valued functions on $\mathbb{R}^{d}$,
the $W_{j}(t)$, $j=1,2,\ldots,m$ are
independent Wiener processes
on a probability space $(\Omega,\cal{F},\mathbb{P})$ with a filtration
$({\cal F}_{t})_{t\geq 0}$
and $\bmath{y}_{0}$ is independent of $W_{j}(t)-W_{j}(0)$.
Especially here, note that $F_{j}$, $j=0,1,\ldots,m$, are
constant matrices and they are noncommutative.

Let $\bmath{y}_{n}$ denote a discrete
approximation to the solution $\bmath{y}(t_{n})$ of \eqref{eq:sde}
for an equidistant grid point $t_{n}\defeq nh$ ($n=1,2,\ldots,M$) with
step size $h=T/M<1$ ($M$ is a positive integer).
In addition, let us introduce the notations
\begin{eqnarray*}
  &&
  \Omega^{[1]}(t_{n},t_{n+1})
  =\left(F_{0}-\frac{1}{2}\sum_{j=1}^{m}F_{j}^{2}
  \right)h
  +\sum_{j=1}^{m}F_{j}\del{W}_{j}(t_{n},t_{n+1}),
  \\
  &&
  \Omega^{[2]}(t_{n},t_{n+1})
  =\Omega^{[1]}(t_{n},t_{n+1})
  +\frac{1}{2}\sum_{i=1}^{m}\sum_{j=i+1}^{m}
  (F_{i}F_{j}-F_{j}F_{i})
  \left(I_{j,i}(t_{n},t_{n+1})-I_{i,j}(t_{n},t_{n+1})
  \right).
\end{eqnarray*}
Then, the Magnus-type Euler and Milstein methods are
respectively given as follows
\cite{Yang:2021}:
\begin{eqnarray}
  &&
  \bmath{y}_{n+1}=\exp\left(\Omega^{[1]}(t_{n},t_{n+1})
  \right)
  \left\{
  \bmath{y}_{n}+\tilde{\bmath{g}}_{0}(\bmath{y}_{n})h
  +\sum_{j=1}^{m}\bmath{g}_{j}(\bmath{y}_{n})
  \del{W}_{j}(t_{n},t_{n+1})
  \right\},
  \label{eq:ME}
  \\
  &&
  \bmath{y}_{n+1}=\exp\left(\Omega^{[2]}(t_{n},t_{n+1})
  \right)
  \left\{
  \bmath{y}_{n}+\tilde{\bmath{g}}_{0}(\bmath{y}_{n})h
  +\sum_{j=1}^{m}\bmath{g}_{j}(\bmath{y}_{n})
  \del{W}_{j}(t_{n},t_{n+1})
  \right.
  \nonumber
  \\
  &&
  \makebox[13.5em]{}
  \left.
  +\sum_{i,j=1}^{m}\bmath{H}_{ij}(\bmath{y}_{n})I_{i,j}(t_{n},t_{n+1})
  \right\},
  \label{eq:MM}
\end{eqnarray}
where
$
\tilde{\bmath{g}}_{0}(\bmath{y})
=\bmath{g}_{0}(\bmath{y})-\sum_{j=1}^{m}F_{j}\bmath{g}_{j}(\bmath{y})
$ and
$\bmath{H}_{i,j}(\bmath{y})
=\bmath{g}_{i}^{\prime}(\bmath{y})
\bigl(F_{j}\bmath{y}+\bmath{g}_{j}(\bmath{y})
\bigr)-F_{j}\bmath{g}_{j}(\bmath{y})
$.

For the linear stability analysis of our methods,
suppose that $\bmath{g}_{j}(\bmath{y})=\bmath{0}$,
$j=0,1,\ldots,m$, in \eqref{eq:sde}.
If we set
\begin{equation}
  \qquad F_{0}=
  \begin{bmatrix}
      \lambda & 0 \\
      0       & \lambda
  \end{bmatrix},
  \qquad
  F_{1}=
  \begin{bmatrix}
      \sigma_{1} & 0 \\
      0 & -\sigma_{1}
  \end{bmatrix},
  \qquad
  F_{2}=
  \begin{bmatrix}
      0 & \sigma_{2} \\
      \sigma_{2} & 0
  \end{bmatrix}
  \label{eq:matrix_non_commutative_test_SDE}
\end{equation}
for $d=m=2$ and nonzero real numbers $\lambda$, $\sigma_{1}$
and $\sigma_{2}$,
then we have the noncommutative test SDE that
Buckwar and Sickenberger \cite{Buckwar:2012} proposed:
\begin{equation}
  \d\bmath{y}(t)=
  F_{0}\bmath{y}(t)\d t
  +F_{1}\bmath{y}(t)\d W_{1}(t)
  +F_{2}\bmath{y}(t)\d W_{2}(t),
  \quad \bmath{y}(0)=\bmath{y}_{0},
  \label{eq:non_commutative_test_SDE}
\end{equation}
in which the zero solution is asymptotically mean square (MS) stable
if and only if
\begin{equation}
  2\lambda+\sigma_{1}^{2}+\sigma_{2}^{2}<0.
  \label{ieq:MS_stab}
\end{equation}
We use the test SDE in our analysis.

In the rest of this subsection, for simplicity we will use the notations
$\del{W}_{j}=\del{W}_{j}(t_{n},t_{n+1})$ and
$A_{i,j}=A_{i,j}(t_{n},t_{n+1})=I_{i,j}(t_{n},t_{n+1})-I_{j,i}(t_{n},t_{n+1})$
without indicating the dependence of $t_{n}$ and $t_{n+1}$,
if it is obvious from the context.
From \eqref{eq:matrix_non_commutative_test_SDE},
\[
\Omega^{[1]}(t_{n},t_{n+1})=\tilde{F}_{0}h+F_{1}\del{W}_{1}+F_{2}\del{W}_{2},
\qquad
\Omega^{[2]}(t_{n},t_{n+1})=\Omega^{[1]}(t_{n},t_{n+1})+GA_{1,2},
\]
where
\[
\tilde{F}_{0}=
\begin{bmatrix}
  \lambda-\frac{1}{2}(\sigma_{1}^{2}+\sigma_{2}^{2}) & 0
  \\
  0 & \lambda-\frac{1}{2}(\sigma_{1}^{2}+\sigma_{2}^{2})
\end{bmatrix},
\qquad
G=
\begin{bmatrix}
    0 & -\sigma_{1}\sigma_{2}
    \\
    \sigma_{1}\sigma_{2} & 0
\end{bmatrix}.
\]
Thus, \eqref{eq:ME} and \eqref{eq:MM} give
the amplification factors
\begin{eqnarray*}
  &&
  R_{E}=\exp(\tilde{F}_{0}h+F_{1}\del{W}_{1}+F_{2}\del{W}_{2}),
  \\
  &&
  R_{M}=\exp(\tilde{F}_{0}h+F_{1}\del{W}_{1}+F_{2}\del{W}_{2}+GA_{1,2}),
\end{eqnarray*}
respectively, as the nonlinear functions
$\bmath{g}_{j}(\bmath{y})$, $j=0,1,\ldots,m$, vanish.

First, let us start with the linear stability analysis
for the Magnus-type Euler method.
We have the following lemma.

\begin{lemma}
  \label{lem:stab_mat_ME}
  When \eqref{eq:ME} is applied to \eqref{eq:non_commutative_test_SDE},
  the stability matrix is expressed as
  \begin{equation}
    E\left[R_{E}^{\top}R_{E}
      \right]
    =\e^{2p-q_{1}-q_{2}}
    \left(
    \sum_{n=0}^{\infty}\varphi_{E}(n)
    \right)I_{d},
    \label{eq:stab_matrix_ME}
  \end{equation}
where $p=\lambda h$, $q_{i}=\sigma_{i}^{2}h$ $(i=1,2)$
and $I_{d}$ stands for the identity matrix,
and where
\begin{equation}
  \varphi_{E}(n)
  =2^{n}
  \sum_{k=0}^{n}
  \frac{{}_{n}C_{k}}{{}_{2n}C_{2k}}
  \frac{q_{1}^{n-k}q_{2}^{k}}{(n-k)!(k!)}.
  \label{eq:diag_component_ME}
\end{equation}
\end{lemma}
\begin{proof}
Noting that
$\tilde{F}_{0}^{\top}=\tilde{F}_{0},F_{1}^{\top}=F_{1}$ and $F_{2}^{\top}=F_{2}$
as well as $\tilde{F}_{0}h$
and $F_{1}\del{W}_{1}+F_{2}\del{W}_{2}$ are commutative,
we have
\[
R_{E}^{\top}R_{E}=\exp(2\tilde{F}_{0}h)
\exp(2(F_{1}\del{W}_{1}+F_{2}\del{W}_{2})).
\]
As $(2(F_{1}\del{W}_{1}+F_{2}\del{W}_{2}))^{2n}
=2^{2n}(\sigma_{1}^{2}\del{W}_{1}^{2}+\sigma_{2}^{2}\del{W}_{2}^{2})^{n}I_{d}
$
for any nonnegative integer $n$, we obtain
\[
E\left[(2(F_{1}\del{W}_{1}+F_{2}\del{W}_{2}))^{2n}
  \right]
=2^{2n}E\left[(\sigma_{1}^{2}\del{W}_{1}^{2}+\sigma_{2}^{2}\del{W}_{2}^{2})^{n}
  \right]I_{d},
\]
whereas $E[(2(F_{1}\del{W}_{1}+F_{2}\del{W}_{2}))^{2n+1}]=O$.
Here,
\begin{eqnarray*}
  &&
  E\left[(\sigma_{1}^{2}\del{W}_{1}^{2}+\sigma_{2}^{2}\del{W}_{2}^{2})^{n}
    \right]
  =
  \sum_{k=0}^{n}{}_{n}C_{n-k}\sigma_{1}^{2(n-k)}
  E\left[\del{W}_{1}^{2(n-k)}
    \right]\sigma_{2}^{2k}
  E\left[\del{W}_{2}^{2k}
    \right]
  \\
  &&
  \makebox[11em]{}
  =\frac{n!}{2^{n}}
  \sum_{k=0}^{n}{}_{2(n-k)}C_{n-k}q_{1}^{n-k}({}_{2k}C_{k})q_{2}^{k}.
\end{eqnarray*}
From these,
\[
E\left[\exp(2(F_{1}\del{W}_{1}+F_{2}\del{W}_{2}))
  \right]
=\sum_{n=0}^{\infty}\frac{1}{(2n)!}
E\left[(2(F_{1}\del{W}_{1}+F_{2}\del{W}_{2}))^{2n}
  \right]
=\varphi_{E}(n).
\]
%
\end{proof}

Let us define A-stability of numerical
methods for the noncommutative
SDE \eqref{eq:non_commutative_test_SDE}
with \eqref{eq:matrix_non_commutative_test_SDE}
and nonzero $\lambda,\sigma_{1},\sigma_{2}$.

\begin{definition}
  A numerical method is said to be A-stable in the MS for
  \eqref{eq:non_commutative_test_SDE}
  if $E[\bmath{y}_{n}^{\top}\bmath{y}_{n}]\to 0$ as $n\to\infty$
  for any step size $h$ whenever the nonzero parameters
  $\lambda,\sigma_{1},\sigma_{2}$ satisfy \eqref{ieq:MS_stab}
  and $E[\bmath{y}_{0}^{\top}\bmath{y}_{0}]<\infty$,
  where $\bmath{y}_{n}$ is a numerical solution given by
  the method when applied to \eqref{eq:non_commutative_test_SDE}.
\end{definition}

As \eqref{eq:diag_component_ME} is a symmetric polynomial in $q_{1},q_{2}$,
let us suppose that $q_{1}\geq q_{2}>0$
without loss of generality and denote $q_{2}/q_{1}$ by $x$.
Then, \eqref{eq:diag_component_ME} is expressed as
\[
\varphi_{E}(n)
=(2q_{1})^{n}
\sum_{k=0}^{n}
\frac{{}_{n}C_{k}}{{}_{2n}C_{2k}}
\frac{x^{k}}{(n-k)!(k!)},
\]
whereas \eqref{ieq:MS_stab} is expressed as
$2p+q_{1}(1+x)<0$.
Regarding the stability of \eqref{eq:ME}, we have the following theorem.

\begin{theorem}
  \label{th:ME}
  The Magnus-type Euler method is A-stable
  in the MS for \eqref{eq:non_commutative_test_SDE}.
\end{theorem}
\begin{proof}
As $\del{W}_{1}(t_{n},t_{n+1})$ and $\del{W}_{2}(t_{n},t_{n+1})$
are independent of ${\cal F}_{t_{n}}$, we have
\[
E[\bmath{y}_{n+1}^{\top}\bmath{y}_{n+1}]
=E[\bmath{y}_{n}^{\top}E[R_{E}^{\top}R_{E}\mid{\cal F}_{t_{n}}]\bmath{y}_{n}]
=E[\bmath{y}_{n}^{\top}E[R_{E}^{\top}R_{E}]\bmath{y}_{n}].
\]
From this and \eqref{eq:stab_matrix_ME},
\begin{equation}
  E[\bmath{y}_{n+1}^{\top}\bmath{y}_{n+1}]
  =\e^{2p-q_{1}(1+x)}
  \left(
  \sum_{n=0}^{\infty}\varphi_{E}(n)
  \right)
  E[\bmath{y}_{n}^{\top}\bmath{y}_{n}].
  \label{eq:admissible_euler}
\end{equation}
Here, this equation is admissible for any $p,q_{1},x$
since $\e^{2p-q_{1}(1+x)}(\sum_{n=0}^{\infty}\varphi_{E}(n))\geq 0$
holds for them clearly.
As ${}_{n}C_{k}\leq{}_{2n}C_{2k}$, $\varphi_{E}(n)$ is
bounded by $(2q_{1})^{n}V_{n}(x)$,
where
\[
V_{n}(x)=\sum_{l=0}^{n}b(n,l)x^{l},
\qquad b(n,l)=\frac{1}{(n-l)!(l!)}.
\]
Thus,
since
\[
\sum_{n=0}^{\infty}\varphi_{E}(n)
\leq
\sum_{n=0}^{\infty}(2q_{1})^{n}V_{n}(x)
=
\sum_{n=0}^{\infty}
\frac{(2q_{1})^{n}}{n!}(1+x)^{n}
=\e^{2q_{1}(1+x)},
\]
we have
\[
\e^{2p-q_{1}(1+x)}
\left(
\sum_{n=0}^{\infty}\varphi_{E}(n)
\right)
\leq
\e^{2p+q_{1}(1+x)}<1
\]
for any $h>0$ in \eqref{eq:stab_matrix_ME} due to
$2p+q_{1}(1+x)<0$ from \eqref{ieq:MS_stab}.
%
\end{proof}

\begin{remark}
  The authors in \cite{Buckwar:2012} have chosen
  another approach on their stability analysis.
  They use the following form:
  \[
  E\left[
  {\rm vec}\left(\bmath{y}_{n+1}\bmath{y}_{n+1}^{\top}
  \right)
  \right]
  =E\left[
  R_{E}\otimes R_{E}
  \right]
  E\left[
    {\rm vec}\left(\bmath{y}_{n}\bmath{y}_{n}^{\top}
    \right)
  \right],
  \]
\end{remark}
where ${\rm vec}(A)$ denotes a vectorisation of a matrix $A$
and $A\otimes B$ denotes the Kronecker product of matrices
$A$ and $B$, and investigate the spectral radius
of $E[R_{E}\otimes R_{E}]$.
On the other hand, in our proof we have avoided dealing
with $E[R_{E}\otimes R_{E}]$, and have mentioned that
\eqref{eq:admissible_euler} is admissible.

\begin{remark}
  In the theorem, we have dealt with $x\in(0,1]$
  for noncommutative noise.
  If $x=0$, then \eqref{eq:non_commutative_test_SDE}
  reduces to an SDE with one scalar noise.
  The Magnus-type Euler method is clearly A-stable also for the
  SDE with one scalar noise.
\end{remark}

Next, let us analyse stability properties
for the Magnus-type Milstein method.
We have the following lemma.

\begin{lemma}
  \label{lem:stab_mat_MM}
  When \eqref{eq:MM} is applied to \eqref{eq:non_commutative_test_SDE},
  the stability matrix is expressed as
  \[
    E\left[R_{M}^{\top}R_{M}
      \right]
    =\e^{2p-q_{1}-q_{2}}
    \left(
    1+\sum_{n=1}^{\infty}\varphi_{M}(n)
    \right)I_{d},
    \]
where
\begin{equation}
  \varphi_{M}(n)
  =\sum_{k=0}^{\tilde{n}}
  (-1)^{k}
  \frac{2^{2(n-k)}}{(2(n-k))!}
  \frac{n-2k}{n-k}
  {}_{n-k}C_{k}
  \sum_{l=k}^{n-k}
      {}_{n-2k}C_{l-k}
      q_{1}^{n-l}
      q_{2}^{l}
      \gamma_{n-k-l,k,l-k}
      \label{eq:diag_component}
\end{equation}
for which $\tilde{n}$ is given by
\begin{equation}
  \tilde{n}=
  \left\{
  \begin{array}{ll}
    r-1
    &
    \quad (n=2r),
    \\
    r
    &
    \quad (n=2r+1).
  \end{array}
  \right.
  \label{eq:k_and_n}
\end{equation}
Here, $r$ stands for an integer.
\end{lemma}
For the proof of this lemma, we refer the reader to
\cref{app:proof_lemma_stab_mat_MM}.

Due to
$\gamma_{n-k-l,k,l-k}=\gamma_{l-k,k,n-k-l}$,
\eqref{eq:diag_component} is a symmetric polynomial in $q_{1},q_{2}$.
By
using \eqref{eq:expect_w1w2DelI21_for_simp}
and \eqref{eq:main_th_rec},
we can rewrite \eqref{eq:diag_component} as
$\varphi_{M}(n)=(2q_{1})^{n}U_{n}(x)$, where
\[
U_{n}(x)
=\sum_{k=0}^{\tilde{n}}
(-1)^{k}
\frac{n-2k}{n-k}
\frac{{}_{n-k}C_{k}}{{}_{2(n-k)}C_{2k}}
\frac{s_{n-2k,k}}{(n-2k)!}
\sum_{l=k}^{n-k}
\frac{({}_{n-2k}C_{l-k})^{2}}{{}_{2(n-2k)}C_{2(l-k)}}
x^{l}.
\]
In a similar way to \eqref{eq:admissible_euler},
we have
\begin{equation}
  E[\bmath{y}_{n+1}^{\top}\bmath{y}_{n+1}]
  =\e^{2p-q_{1}(1+x)}
  \left(
  1+\sum_{n=1}^{\infty}(2q_{1})^{n}U_{n}(x)
  \right)
  E[\bmath{y}_{n}^{\top}\bmath{y}_{n}].
  \label{eq:ms_milstein}
\end{equation}

\begin{remark}
  \label{rem:un}
  For a fixed $l$,
  the monomials with respect to $x^{l}$ and $x^{n-l}$
  appear only for $k=0,1,\ldots,\min(l,\tilde{n})$ in $U_{n}(x)$.
  Noting this, we can express $U_{n}(x)$ in a usual form with coefficients
  $a(n,l),\ l=0,1,\ldots,n$:
  $
    U_{n}(x)=\sum_{l=0}^{n}a(n,l)x^{l},
    $
  where
  \begin{eqnarray*}
    &&
    a(n,l)=
    \sum_{k=0}^{l}\delta_{n,k,l}
    \ \ (l=0,1,\ldots,\tilde{n}),
    \quad
    a(n,l)=
    \sum_{k=0}^{n-l}\delta_{n,k,l}
    \ \ (l=\tilde{n}+2,\tilde{n}+3,\ldots,n),
    \\
    &&
    a(n,\tilde{n}+1)=
    \sum_{k=0}^{\tilde{n}}\delta_{n,k,\tilde{n}+1},
    \quad
    \delta_{n,k,l}
    =
    (-1)^{k}
    \frac{n-2k}{n-k}
    \frac{{}_{n-k}C_{k}}{{}_{2(n-k)}C_{2k}}
    \frac{s_{n-2k,k}}{(n-2k)!}
    \frac{({}_{n-2k}C_{l-k})^{2}}{{}_{2(n-2k)}C_{2(l-k)}}.
  \end{eqnarray*}
  If the coefficients are given,
  the usual form with Horner's rule is helpful for
  the faster calculation of the value of $U_{n}(x)$ at each given $x$.
\end{remark}


For the stability analysis of \eqref{eq:MM}
when it is applied to \eqref{eq:non_commutative_test_SDE},
we investigate $U_{n}(x)$ in detail.

\begin{figure}[t]
  \unitlength 1cm
  \begin{center}
    \begin{picture}(13,3.7)
      %
      \put(0,0){\includegraphics[width=5.8cm,keepaspectratio]{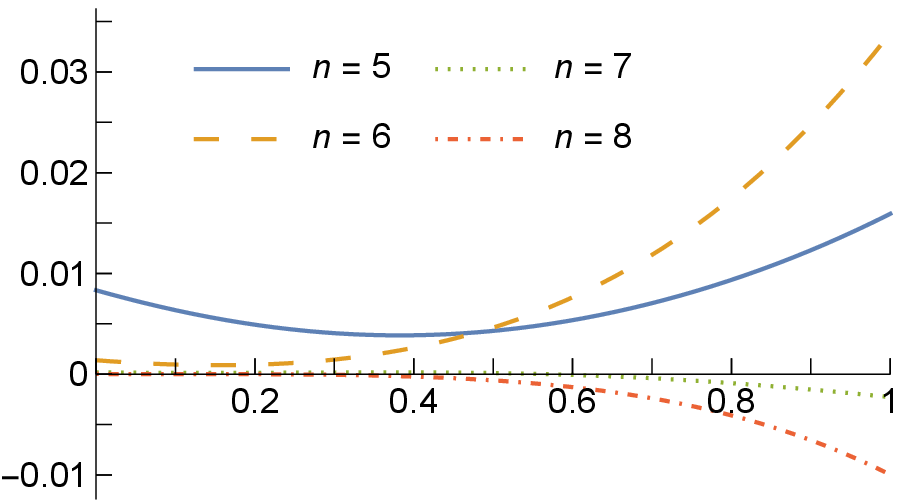}}
      \put(6.7,0){\includegraphics[width=5.8cm,keepaspectratio]{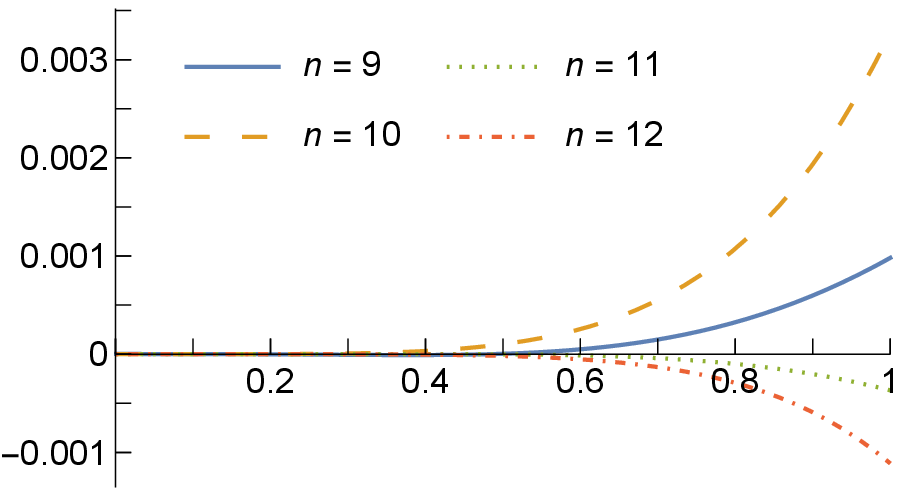}}
      \put(0.2, 3.4){$U_{n}(x)$}
      \put(6.1, 0.8){$x$}
      \put(7.1, 3.4){$U_{n}(x)$}
      \put(12.8, 0.8){$x$}
    \end{picture}
  \end{center}
  \caption{Plots of $U_{n}(x)$ for $n=5,6,\ldots,12$}
  \label{fig:Un}
\end{figure}

\subsubsection{Analysis for the Magnus-type Milstein method}

\begin{figure}[t]
  \unitlength 1cm
  \begin{center}
    \begin{picture}(13.0,3.9)
      %
      \put(0,0){\includegraphics[width=5.8cm,keepaspectratio]{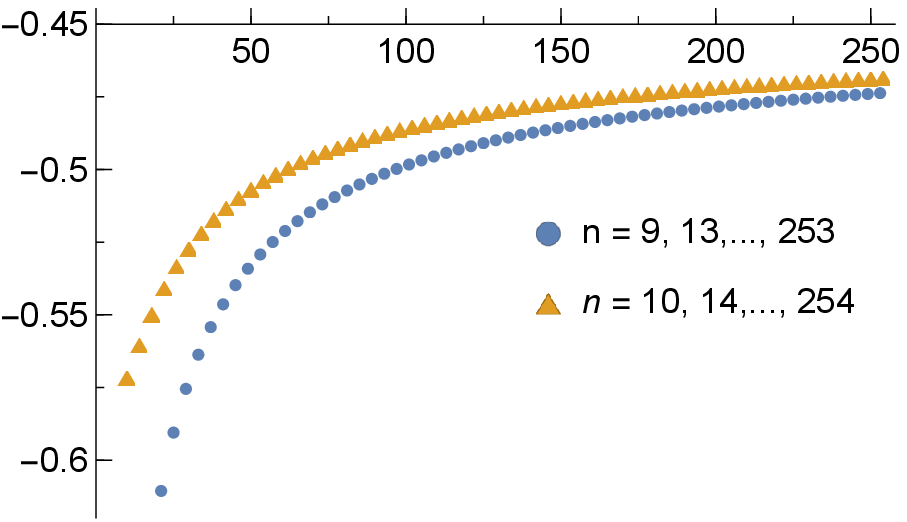}}
      \put(6.8,0){\includegraphics[width=5.8cm,keepaspectratio]{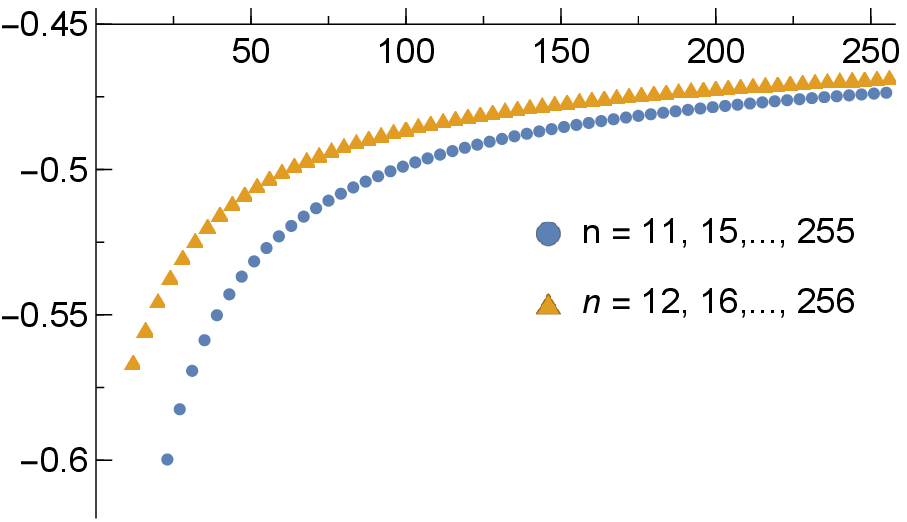}}
      \put(0, 3.6){$\frac{1}{n}\ln U_{n}(1)$}
      \put(6.1, 3.1){$n$}
      \put(6.9, 3.6){$\frac{1}{n}\ln (-U_{n}(1))$}
      \put(12.85, 3.1){$n$}
    \end{picture}
  \end{center}
  \caption{$\frac{1}{n}\ln|U_{n}(1)|$ versus $n$}
  \label{fig:LyapunovExp}
\end{figure}

Let us give the first four examples of $U_{n}(x)$
\begin{eqnarray*}
  &&
  U_{1}(x)=1+x,
  \qquad
  U_{2}(x)=\frac{1}{2}+\frac{1}{3}x+\frac{1}{2}x^{2},
  \qquad
  U_{3}(x)=\frac{1}{6}-\frac{17}{180}x-\frac{17}{180}x^{2}+\frac{1}{6}x^{3},
  \\
  &&
  U_{4}(x)=\frac{1}{24}-\frac{17}{315}x-\frac{23}{756}x^{2}
  -\frac{17}{315}x^{3}+\frac{1}{24}x^{4}
\end{eqnarray*}
and show others in \cref{fig:Un}. In addition, let us
investigate $(1/n)\ln|U_{n}(1)|$
separately for $n=4k+i$, where $i=1,2,3,4$ and
$k=2,3,\ldots,63$,
and show them in \cref{fig:LyapunovExp}.
More concretely, for example,
$(1/n)\ln(U_{n}(1))=-0.473732$ and
$(1/n)\ln(-U_{n}(1))=-0.468831$
for $n=253,256$, respectively.
These lead to $U_{253}(1)=(0.622674)^{253}$
and $U_{256}(1)=-(0.625734)^{256}$
as $\exp(-0.473732)=0.622674$ and $\exp(-0.468831)=0.625734$.
Taking these into account, when $x=1$, 
we can see that the series in \eqref{eq:ms_milstein}
does not even converge for large $q_{1}$ such as $q_{1}=1$.
That is, \eqref{eq:ms_milstein} is inadmissible
for large $q_{1}$ when $x=1$.

\begin{remark}
  If $x=0$, then \eqref{eq:non_commutative_test_SDE}
  reduces to an SDE with one scalar noise.
  The Magnus-type Milstein method is clearly A-stable for the
  SDE with one scalar noise.
\end{remark}

\begin{figure}[t]
\unitlength 1cm
\begin{center}
\begin{picture}(13.0,4.5)
%
\put(0.5,0.8){\includegraphics[width=5.8cm,keepaspectratio]{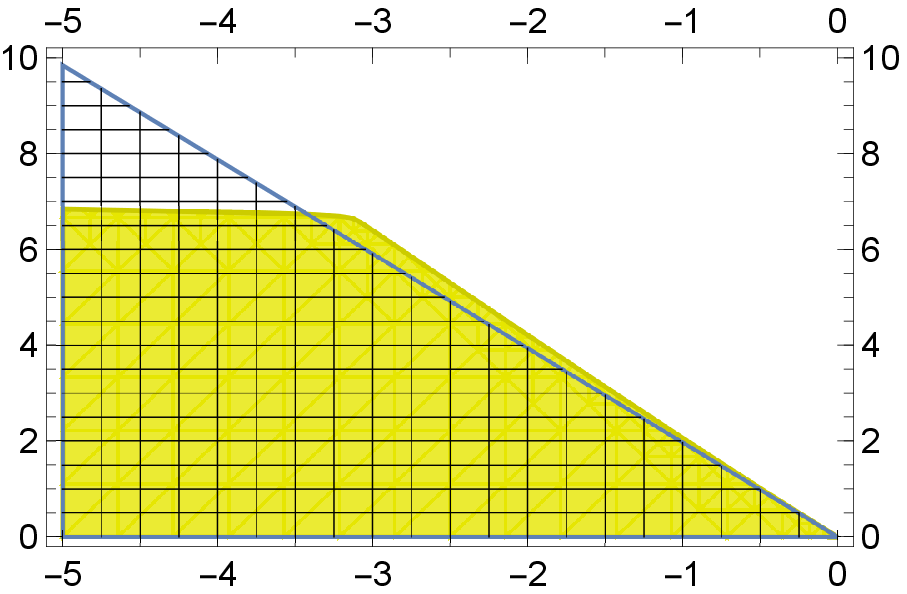}}
\put(7.2,0.8){\includegraphics[width=5.8cm,keepaspectratio]{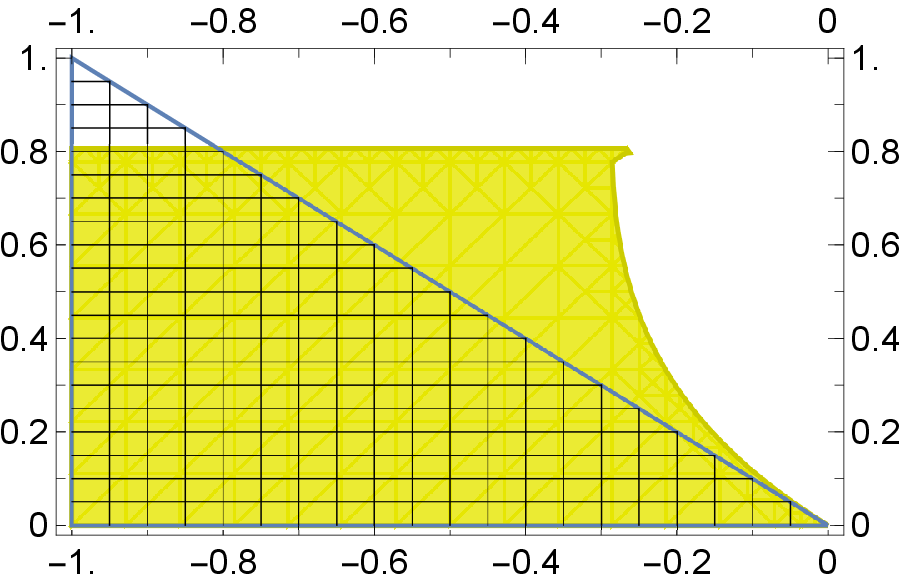}}
\put(0, 3.8){$q_{1}$}
\put(5.3, 0.5){$p$}
\put(2.7,0.0){$x=1/64$}
\put(6.7, 3.8){$q_{1}$}
\put(12.0, 0.5){$p$}
\put(9.6,0.0){$x=1$}
\end{picture}
\caption{Stability regions for the Magnus-type Milstein method
when $x=1/64, 1$}
\label{fig:stab_region}
\end{center}
\end{figure}

Now, let us plot stability regions for the Magnus-type
Milstein method by
\begin{equation}
  0\leq
  \e^{2p-q_{1}(1+x)}
  \left(
  1+\sum_{n=1}^{\infty}(2q_{1})^{n}U_{n}(x)
  \right)
  <1.
  \label{ieq:stab_cond_Mil}
\end{equation}
We utilize a truncated expansion for
the series in the inequality.
\cref{fig:stab_region} indicates
the stability regions by colored part, when
$U_{n}(x)$, $n=1,2,\ldots,256$, are used with $x=1/64,1$.
The other parts enclosed by the mesh indicate
the region that satisfies $2p+q_{1}(1+x)<0$ for $x=1/64,1$.

\begin{remark}
  The results obtained so far allow us to infer the following:
  there might exist a positive real number $q_{0}$ such that
  \eqref{ieq:stab_cond_Mil} holds for any $x\in[0,1]$
  and $0< q_{1}\leq q_{0}$ under the condition
  $2p+q_{1}(1+x)<0$,
  but $q_{0}$ would not be very large.
  %
\end{remark}

\begin{table}[t]
  \footnotesize
  \caption{Logarithm of the error of $\bmath{y}_{n}$ at $t_{n}=1$ to base 2}
  \label{tab:errors}
  \begin{center}
    \begin{tabular}{c|cccccc}
      $\log_{2}h$ & $-1$ & $-2$ & $-3$ & $-4$ & $-5$ & $-6$
      \\
      \hline
      Magnus-type Euler
      & $-7.49$ & $-7.91$ & $-8.42$ & $-8.94$ & $-9.39$ & $-9.88$
      \\
      Magnus-type Milstein
      & $-8.80$ & $-9.78$ & $-10.7$ & $-11.7$ & $-12.7$ & $-13.6$
    \end{tabular}
  \end{center}
\end{table}

\begin{table}[t]
  \footnotesize
  \caption{Arithmetic mean and standard deviation
    of $\|\bmath{y}_{n}\|^{2}$ at $t_{n}=5$
  by the Magnus-type Milstein method when $q_{1}=1/2$ and $x=1$}
  \label{tab:2nd_moment}
  \begin{center}
    \begin{tabular}{c|cccc}
      $p$ & $-0.1$ & $-0.2$ & $-0.4$ & $-0.6$
      \\
      \hline
      $\langle\|\bmath{y}_{n}\|^{2}\rangle$
      & $5.17\times 10$ & $6.99$ & $1.28\times 10^{-1}$ & $2.35\times 10^{-3}$
      \\
      {\rm SD}
      & $2.64\times 10$ & $3.57$ & $6.55\times 10^{-2}$ & $1.20\times 10^{-3}$
    \end{tabular}
  \end{center}
\end{table}

If $q_{1}$ is not very large, it is not difficult to confirm
our theoretical results numerically
(see also
\url{https://github.com/yosh-komori/mag-simul}).
In what follows, let us use a fixed initial
condition $\bmath{y}(0)=[1\ 1]^{\top}$.
When the Magnus-type
methods are applied to \eqref{eq:non_commutative_test_SDE}
with $\lambda=-1/4,\sigma_{1}=1/2,\sigma_{2}=2/5$
and $1000$ independent trajectories are simulated for
a given $h$, \cref{tab:errors} shows errors
$\langle\|\bmath{y}_{n}-\bmath{y}_{n}^{\rm ref}\|^{2}\rangle$ at $t_{n}=1$,
which stands for the arithmetic mean of
$\|\bmath{y}_{n}-\bmath{y}_{n}^{\rm ref}\|^{2}$, and
where $\bmath{y}_{n}^{\rm ref}$ is a reference solution computed
by the Milstein method with $h=2^{-9}$.
From the table, we can see that the methods
achieve the theoretical order of convergence.
When the Magnus-type Milstein method
with $h=2^{-1}$ is applied to \eqref{eq:non_commutative_test_SDE}
with $\lambda=-1/5,-2/5,-4/5,-6/5,\sigma_{1}=1$
and $\sigma_{2}=1$, and $10$ batches of the $10^{5}$ independent trajectories
are simulated for each value of $\lambda$, \cref{tab:2nd_moment}
shows the arithmetic mean and standard deviation (SD)
of $\|\bmath{y}_{n}\|^{2}$ at $t_{n}=5$.
From the table, we can see that the method
numerically reproduces the stability properties
indicated by the right-hand plot in \cref{fig:stab_region}.

On the other hand, as $q_{1}$ becomes larger, situations change.
When we replace $\sigma_{1}=\sigma_{2}=1$
with $\sigma_{1}=\sigma_{2}=3$ and set $\lambda=-0.01$
in the last case,
the setting makes numerical solutions of
the method unstable,
but the numerical value of the second moment by the method with $h=2^{-1}$
is $3.17\times 10^{-14}$ $(9.32\times 10^{-14})$,
where the value in parentheses indicates SD.
It is much smaller than the values in \cref{tab:2nd_moment}.
One of the reasons for this phenomena can be considered as follows.
The MS instability of the method is determined by rare
exploding trajectories. The standard Monte
Carlo approach is bound to miss those rare events
\cite{Ableidinger:2017},
and
$I_{i,j}(t_{n},t_{n+1})$
is approximated
\cite{Wiktorsson:2001}.
The authors in \cite{Ableidinger:2017} have tackled such a problem
in a multi-dimensional linear SDE with diagonal diffusion
matrices, and they have pointed out that
noncommutative cases are more challenging.
We leave this problem as a future work.

\subsection{Mixed moments of Wiener increments and
  stochastic double integrals}

In \cref{subsec:preliminary}, we mentioned
that mixed moments of
$\del{W}_{1}(t,t+1),\del{W}_{2}(t,t+1),I_{1,2}(t,t+1)$
and $I_{2,1}(t,t+1)$ are
obtained from the properties of stochastic
integrals, but it is not easy to seek them when an exponent is large
number.
As another application of our main theorem,
we give such mixed moments.
In the rest of this subsection, for simplicity we will use the notations
$\del{W}_{j}=\del{W}_{j}(t,t+1)$,
$A_{i,j}=A_{i,j}(t,t+1)$ and $I_{i,j}=I_{i,j}(t,t+1)$.

Noting $A_{2,1}=-A_{1,2}$, we have
$I_{1,2}+I_{2,1}=\del{W}_{1}\del{W}_{2}$.
This and $I_{1,2}-I_{2,1}=A_{1,2}$ give
\[
A_{1,2}^{2}
=(\del{W}_{1}\del{W}_{2})^{2}-4I_{1,2}I_{2,1}.
\]
Thus, we have
$
E\left[
  \del{W}_{1}^{2}
  A_{1,2}^{2}
  \right]
=E\left[
  \del{W}_{1}^{4}
  \del{W}_{2}^{2}
  \right]
-4E\left[
  \del{W}_{1}^{2}
  I_{1,2}I_{2,1}
  \right]
$.
As our main theorem gives
\[
E\left[
  \del{W}_{1}^{2}
  A_{1,2}^{2}
  \right]
=\gamma_{1,1,0}
=\frac{5}{3},
\]
we have
$
E\left[
  \del{W}_{1}^{2}
  I_{1,2}I_{2,1}
  \right]
=\frac{1}{3}
$.
Note that this expectation is replaced with $h^{3}/3$
if we change $t+1$ to $t+h$ in $\del{W}_{1},I_{1,2}$ and $I_{2,1}$
\cite[p. 225]{Kloeden:1999}.

Similarly, we have
\begin{eqnarray*}
  &&
  E\left[
    \del{W}_{1}^{4}
    A_{1,2}^{2}
    \right]
  =E\left[
    \del{W}_{1}^{6}
    \del{W}_{2}^{2}
    \right]
  -4E\left[
    \del{W}_{1}^{4}
    I_{1,2}I_{2,1}
    \right],
  \\
  &&
  E\left[
    \del{W}_{1}^{2}
    A_{1,2}^{2}
    \del{W}_{2}^{2}
    \right]
  =E\left[
    \del{W}_{1}^{4}
    \del{W}_{2}^{4}
    \right]
  -4E\left[
    \del{W}_{1}^{2}
    \del{W}_{2}^{2}
    I_{1,2}I_{2,1}
    \right].
\end{eqnarray*}
On the other hand, our main theorem
and \eqref{eq:expect_w1w2DelI21_for_simp} give
\[
E\left[
  \del{W}_{1}^{4}
  A_{1,2}^{2}
  \right]
=\gamma_{2,1,0}
=7,
\qquad
E\left[
  \del{W}_{1}^{2}
  A_{1,2}^{2}
  \del{W}_{2}^{2}
  \right]
=\gamma_{1,1,1}
=\frac{{}_{2}C_{1}}{{}_{4}C_{2}}\gamma_{2,1,0}
=\frac{7}{3}.
\]
Thus, we have
$
E\left[
  \del{W}_{1}^{4}
  I_{1,2}I_{2,1}
  \right]
=2
$ and
$
\left[
  \del{W}_{1}^{2}
  \del{W}_{2}^{2}
  I_{1,2}I_{2,1}
  \right]
=5/3
$.
In a similar way, we can easily obtain
the following mixed moments and others.
\begin{eqnarray*}
  &&
  E\left[
    \del{W}_{1}^{6}
    I_{1,2}I_{2,1}
    \right]
  =15,
  \quad
  E\left[
    \del{W}_{1}^{4}
    \del{W}_{2}^{2}
    I_{1,2}I_{2,1}
    \right]
  =9,
  \quad
  E\left[
    \del{W}_{1}^{2}
    I_{1,2}^{2}I_{2,1}^{2}
    \right]
  =\frac{49}{20},
  \\
  &&
  E\left[
    \del{W}_{1}^{8}
    I_{1,2}I_{2,1}
    \right]
  =140,
  \quad
  E\left[
    \del{W}_{1}^{6}
    \del{W}_{2}^{2}
    I_{1,2}I_{2,1}
    \right]
  =65,
  \quad
  E\left[
    \del{W}_{1}^{4}
    \del{W}_{2}^{4}
    I_{1,2}I_{2,1}
    \right]
  =48,
  \\
  &&
  E\left[
    \del{W}_{1}^{4}
    I_{1,2}^{2}I_{2,1}^{2}
    \right]
  =\frac{339}{20},
  \quad
  E\left[
    \del{W}_{1}^{2}
    \del{W}_{2}^{2}
    I_{1,2}^{2}I_{2,1}^{2}
    \right]
  =\frac{679}{60}.
\end{eqnarray*}
\section{Concluding remarks}
\label{sec:conclusions}
We have derived a way of obtaining the expectation of monomial
$
\gamma_{n,k,l}
=E\left[\{\del{W}_{1}(t,t+1)\}^{2n}
  \{A_{1,2}(t,t+1)\}^{2k}
  \{\del{W}_{2}(t,t+1)\}^{2l}
  \right]
$
for nonnegative integers $n,k$ and $l$.
Due to \eqref{eq:expect_w1w2DelI21_for_simp},
it is essential to obtain $\gamma_{n,k,0}$.
We have derived the three types of formulae which give it:
a recursive formula, an explicit formula,
and a simpler formula that includes
the mixed moment generating function of
$\{\del{W}_{1}(t,t+1)\}^{2}$
and $\{A_{1,2}(t,t+1)\}^{2}$.

If a numerical method contains
$I_{i,j}(t_{n},t_{n+1})=(1/2)\Delta{W}_{i}(t_{n},t_{n+1})\Delta{W}_{j}(t_{n},t_{n+1})
+(1/2)A_{i,j}(t_{n},t_{n+1})$ and
if we use the Taylor expansion of a function $f(y_{n+1})$ centered at $y_{n}$
for the analysis of the method,
then we often need to calculate $\gamma_{n,k,l}$.
One of such methods is the Magnus-type Milstein method.
Utilizing the formulae, we have analysed the stability
properties of the method when it is applied to the
noncommutative test SDE \eqref{eq:non_commutative_test_SDE}.
As a result, our analysis
has
shown that
the method cannot achieve A-stability for the test SDE.

Finally, we make the following remarks.
If $\sigma_{1}h$ or $\sigma_{2}h$ is not very large for the parameters
$\sigma_{1},\sigma_{2}$ in \eqref{eq:non_commutative_test_SDE},
it is not difficult to confirm our theoretical results numerically.
On the other hand, if they are large, situations change.
Due to the term $\exp\left(\Omega^{[2]}(t_{n},t_{n+1})\right)$ in the method,
the MS instability of the method is determined by rare
exploding trajectories. The standard Monte
Carlo approach is bound to miss those rare events.
This is another challenging issue in noncommutative cases.
Thus, we will consider this issue in future work.

\appendix
\section{Proof of Lemma \ref{lem:sum_sHat}}
\label{app:proof_lemma_sum_sHat}


When $L=0$, it is trivial that the lemma holds.
As $\hat{s}_{n,L}(k)=\hat{s}_{n,k}(k)$ for $L>k$,
we can suppose that $L\leq k$ without loss of generality.
Then, similarly to a part of the proof of
Lemma \ref{lem:explicit_formula_rn},
\begin{eqnarray*}
  &&
  \makebox[1em]{}
  \sum_{j=1}^{L}\beta_{n,j}\hat{s}_{n,L}(k-j)
  \\
  &&
  =
  \sum_{\hat{l}_{1}+2l_{2}+\cdots+Ll_{L}=k}
  \left\{
    \hat{l}_{1}
    \frac{\beta_{n,1}^{\hat{l}_{1}}}{\hat{l}_{1}!}
    \prod_{\stackrel{\scriptstyle j=1}{j \neq 1}}^{L}
    \frac{\beta_{n,j}^{l_{j}}}{j^{l_{j}}(l_{j}!)}
    \right\}
  +
  \sum_{l_{1}+2\hat{l}_{2}+\cdots+Ll_{L}=k}
  \left\{
  2
  \hat{l}_{2}
  \frac{\beta_{n,2}^{\hat{l}_{2}}}{
    2^{\hat{l}_{2}}
    (\hat{l}_{2}!)}
  \prod_{\stackrel{\scriptstyle j=1}{j \neq 2}}^{L}
  \frac{\beta_{n,j}^{l_{j}}}{j^{l_{j}}(l_{j}!)}
  \right\}
  \\
  &&
  \makebox[1em]{}
  +\cdots
  +\sum_{l_{1}+2l_{2}+\cdots+(L-1)l_{L-1}+L\hat{l}_{L}=k}
  \left\{
    L\hat{l}_{L}
    \frac{\beta_{n,L}^{\hat{l}_{L}}}{L^{\hat{l}_{L}}(\hat{l}_{L}!)}
    \prod_{\stackrel{\scriptstyle j=1}{j \neq L}}^{L}
    \frac{\beta_{n,j}^{l_{j}}}{j^{l_{j}}(l_{j}!)}
    \right\}
  \\
  &&
  =
  \sum_{l_{1}+2l_{2}+\cdots+Ll_{L}=k}
  \left\{
  \left(l_{1}+2l_{2}+\cdots+Ll_{L}
  \right)\prod_{j=1}^{L}\frac{\beta_{n,j}^{l_{j}}}{j^{l_{j}}(l_{j}!)}
  \right\}
  =k\hat{s}_{n,L}(k),
\end{eqnarray*}
where $\hat{l}_{1},\hat{l}_{2},\ldots,\hat{l}_{L}$
denote positive integers.
This completes the proof.


\section{Proof of Lemma \ref{lem:stab_mat_MM}}
\label{app:proof_lemma_stab_mat_MM}


Noting that
$\tilde{F}_{0}^{\top}=\tilde{F}_{0},F_{1}^{\top}=F_{1},F_{2}^{\top}=F_{2}$
and $G^{\top}=-G$
as well as $\tilde{F}_{0}h$
and $F_{1}\del{W}_{1}+F_{2}\del{W}_{2}+GA_{1,2}$ are commutative,
we have
\[
R_{M}^{\top}R_{M}=\exp(2\tilde{F}_{0}h)\e^{P}\e^{Q},
\]
where
$
P=F_{1}\del{W}_{1}+F_{2}\del{W}_{2}-GA_{1,2}
$ and
$
Q=F_{1}\del{W}_{1}+F_{2}\del{W}_{2}+GA_{1,2}
$.
Although $P$ and $Q$ are noncommutative,
they satisfy the equations
\begin{eqnarray*}
  &&
  PQ=(\sigma_{1}^{2}\del{W}_{1}^{2}+\sigma_{2}^{2}\del{W}_{2}^{2}
  +(\sigma_{1}\sigma_{2})^{2}A_{1,2}^{2})I_{d}
  +2F_{1}GA_{1,2}\del{W}_{1}
  +2F_{2}GA_{1,2}\del{W}_{2}.
  \\
  &&
  P^{2k}
  =Q^{2k}
  =(\sigma_{1}^{2}\del{W}_{1}^{2}+\sigma_{2}^{2}\del{W}_{2}^{2}
  -(\sigma_{1}\sigma_{2})^{2}A_{1,2}^{2})^{k}I_{d}
\end{eqnarray*}
for any positive integer $k$.
From these and the properties of the random variables,
we have
\begin{eqnarray*}
  &&
  E[PQ]=
  E\left[
    (\sigma_{1}^{2}\del{W}_{1}^{2}+\sigma_{2}^{2}\del{W}_{2}^{2}
    +(\sigma_{1}\sigma_{2})^{2}A_{1,2}^{2})
    \right]I_{d}
  =E\left[
    P^{2}+2(\sigma_{1}\sigma_{2})^{2}A_{1,2}^{2}I_{d}
    \right],
  \\
  &&
  E\left[
    P^{2k}
    \right]
  =
  E\left[
    Q^{2k}
    \right]
  =
  E\left[
    (\sigma_{1}^{2}\del{W}_{1}^{2}+\sigma_{2}^{2}\del{W}_{2}^{2}
    -(\sigma_{1}\sigma_{2})^{2}A_{1,2}^{2})^{k}
    \right]I_{d},
  \\
  &&
  E\left[
    P^{2k+1}
    \right]
  =E\left[
    Q^{2k+1}
    \right]
  =O,
  \qquad
  E\left[
    P^{k}Q^{l}
    \right]
  =O\quad({\rm if}\ k+l\ {\rm is\ an\ odd\ number}),
  \\
  &&
  E\left[
    P^{2k}PQQ^{2l}
    \right]
  =
  E\left[
    P^{2k}PQP^{2l}
    \right]
  =
  E\left[
    P^{2(k+l+1)}+2(\sigma_{1}\sigma_{2})^{2}A_{1,2}^{2}P^{2(k+l)}
    \right].
\end{eqnarray*}
The application of these relationships to
\begin{eqnarray*}
  &&
  \makebox[1em]{}
  E\left[
    \e^{P}\e^{Q}
    \right]
  \\
  &&
  =
  E\left[
    I_{d}+\frac{1}{2!}(P^{2}+2PQ+Q^{2})
    +\frac{1}{4!}(P^{4}+4P^{3}Q+6P^{2}Q^{2}+4PQ^{3}+Q^{4})
    \right.
    \\
    &&
    \makebox[3em]{}
    \left.
    +\frac{1}{6!}(P^{6}+6P^{5}Q+15P^{4}Q^{2}+20P^{3}Q^{3}
    +15P^{2}Q^{4}+6PQ^{5}+Q^{6})
    +\cdots
    \right]
\end{eqnarray*}
leads to
\begin{eqnarray*}
  &&
  \makebox[1em]{}
  E\left[
    \e^{P}\e^{Q}
    \right]
  =E\left[
    I_{d}+\frac{2^{2}}{2!}(\sigma_{1}^{2}\del{W}_{1}^{2}
    +\sigma_{2}^{2}\del{W}_{2}^{2})I_{d}
    +\frac{2^{4}}{4!}(\sigma_{1}^{2}\del{W}_{1}^{2}
    +\sigma_{2}^{2}\del{W}_{2}^{2})P^{2}
    \right.
  \\
  &&
  \makebox[8.5em]{}
  \left.
  +\frac{2^{6}}{6!}(\sigma_{1}^{2}\del{W}_{1}^{2}
  +\sigma_{2}^{2}\del{W}_{2}^{2})P^{4}
  +\cdots
  \right].
\end{eqnarray*}
The substitution of
$P^{2}=
\left(
\sigma_{1}^{2}\del{W}_{1}^{2}+\sigma_{2}^{2}\del{W}_{2}^{2}
-(\sigma_{1}\sigma)^{2}A_{1,2}^{2}
\right)I_{d}
$ into this and simplification yield
\[
E\left[
  \e^{P}\e^{Q}
  \right]
=
E\left[
  1+\sum_{n=1}^{\infty}
  \frac{2^{2n}}{(2n)!}
  \sum_{k=0}^{n-1}
      {}_{n-1}C_{k}
      \left(\sigma_{1}^{2}\del{W}_{1}^{2}+\sigma_{2}^{2}\del{W}_{2}^{2}
      \right)^{n-k}
      \left\{-(\sigma_{1}\sigma_{2})^{2}A_{1,2}^{2}
      \right\}^{k}
      \right]I_{d}.
\]
Noting that
\[
E\left[
  \left(\sigma_{1}^{2}\del{W}_{1}^{2}+\sigma_{2}^{2}\del{W}_{2}^{2}
      \right)^{n-k}
      \right]=O(h^{n-k}),
\qquad
E\left[
  \left\{-(\sigma_{1}\sigma_{2})^{2}A_{1,2}^{2}
  \right\}^{k}
  \right]=O(h^{2k}),
\]
we can see that the coefficient of $O(h^{n})$ in
$E\left[\e^{P}\e^{Q}\right]$, say $\varphi_{M}(n)$, is given by
\begin{eqnarray}
  &&
  E\left[
    \frac{2^{2n}}{(2n)!}
         {}_{n-1}C_{0}
         \left(\sigma_{1}^{2}\del{W}_{1}^{2}+\sigma_{2}^{2}\del{W}_{2}^{2}
         \right)^{n}
         \right.
         \nonumber
         \\
         &&
         \makebox[1.5em]{}
         +\frac{2^{2(n-1)}}{(2(n-1))!}
         {}_{n-2}C_{1}
         \left(\sigma_{1}^{2}\del{W}_{1}^{2}+\sigma_{2}^{2}\del{W}_{2}^{2}
         \right)^{n-2}
         \left\{-(\sigma_{1}\sigma_{2})^{2}A_{1,2}^{2}
         \right\}
         \nonumber
         \\
         &&
         \makebox[1.5em]{}
         +\cdots
         \left.
         +\frac{2^{2(n-k)}}{(2(n-k))!}
         {}_{n-k-1}C_{k}
         \left(\sigma_{1}^{2}\del{W}_{1}^{2}+\sigma_{2}^{2}\del{W}_{2}^{2}
         \right)^{n-2k}
         \left\{-(\sigma_{1}\sigma_{2})^{2}A_{1,2}^{2}
         \right\}^{k}
         \right],
  \label{eq:diag_term}
\end{eqnarray}
where $n-2k\geq 1$.
Here, since
\begin{eqnarray*}
  &&
  \makebox[1em]{}
  \left(
  \sigma_{1}^{2}\del{W}_{1}^{2}+\sigma_{2}^{2}\del{W}_{2}^{2}
  \right)^{n-2k}
  \left\{-(\sigma_{1}\sigma_{2})^{2}A_{1,2}^{2}
  \right\}^{k}
  \\
  &&
  =
  (-1)^{k}
  \sum_{l=0}^{n-2k}
      {}_{n-2k}C_{l}
      \sigma_{1}^{2(n-k-l)}
      \sigma_{2}^{2(k+l)}
      \del{W}_{1}^{2(n-2k-l)}
      \del{W}_{2}^{2l}
      A_{1,2}^{2k}
\end{eqnarray*}
and
\[
E\left[
  \left\{
  \del{W}_{1}(t_{n},t_{n+1})
  \right\}^{2(n-2k-l)}
  \left\{
  \del{W}_{2}(t_{n},t_{n+1})
  \right\}^{2l}
  \left\{
  A_{1,2}(t_{n},t_{n+1})
  \right\}^{2k}
  \right]
=\gamma_{n-2k-l,k,l}h^{n},
\]
we have
\begin{eqnarray*}
  &&
  \makebox[1em]{}
  \frac{2^{2(n-k)}}{(2(n-k))!}{}_{n-k-1}C_{k}
  E\left[
    \left(
    \sigma_{1}^{2}\del{W}_{1}^{2}+\sigma_{2}^{2}\del{W}_{2}^{2}
    \right)^{n-2k}
    \left\{-(\sigma_{1}\sigma_{2})^{2}A_{1,2}^{2}
    \right\}^{k}
    \right]
  \\
  &&
  =
  (-1)^{k}
  \frac{2^{2(n-k)}}{(2(n-k))!}
  \frac{n-2k}{n-k}
  {}_{n-k}C_{k}
  \sum_{l=k}^{n-k}
      {}_{n-2k}C_{l-k}
      q_{1}^{n-l}
      q_{2}^{l}
      \gamma_{n-k-l,k,l-k}
\end{eqnarray*}
in \eqref{eq:diag_term}.
Here, remember that $q_{i}=\sigma_{i}^{2}h$ ($i=1,2$).
Thus, $\varphi_{M}(n)$ is rewritten as
\[
\varphi_{M}(n)
=\sum_{k=0}^{\tilde{n}}
(-1)^{k}
\frac{2^{2(n-k)}}{(2(n-k))!}
\frac{n-2k}{n-k}
{}_{n-k}C_{k}
\sum_{l=k}^{n-k}
    {}_{n-2k}C_{l-k}
    q_{1}^{n-l}
    q_{2}^{l}
    \gamma_{n-k-l,k,l-k}.
\]
As $k$ must satisfy $n-2k\geq 1$,
$\tilde{n}$ is given by \eqref{eq:k_and_n}.


%

\section*{Acknowledgments}
The authors would like to thank referees for their
comments which helped to improve the earlier version of
this paper.

\bibliographystyle{siamplain}
%

%
\end{document}